\documentclass[a4paper,11pt]{article}

\usepackage[a4paper,left=3.0cm,right=3.0cm]{geometry}

\usepackage{amsmath}
\usepackage{amssymb}
\usepackage{mathtools}
\usepackage{amsfonts}  
\usepackage{bbm}
\usepackage{mathrsfs} 
\usepackage{xifthen} 
\usepackage[]{xcolor}
\usepackage{bigints}
\usepackage{booktabs}
\usepackage{para list}
\usepackage{stmaryrd}
\usepackage[normalem]{ulem}

\usepackage{graphicx}
\usepackage{newtxtext}
\usepackage{newtxmath}
\usepackage{hyperref}
\hypersetup{
    colorlinks = true,
    urlcolor   = blue,
    citecolor  = black,
}
\newcommand{\RomanNumeralCaps}[1]
\linenumbers

\usepackage{mathrsfs}
\usepackage{scalerel}

\DeclareMathAlphabet{\Mymathbb}{U}{bbold}{m}{n}
\DeclareMathAlphabet{\mathpzc}{OT1}{pzc}{m}{it}

\newtheorem{remark}{Remark}
\newtheorem{problem}{Problem}

\newcommand{\timeSymbol}{t}
\newcommand{\Time}{\timeSymbol} 

\newcommand{\FirstSymbol}{\mathcal{G}}
\newcommand{\First}[1][]{
  \ifthenelse{\equal{#1}{}}
  {\FirstSymbol}
  {\FirstSymbol_{\scriptscriptstyle{#1}}}
}

\newcommand{\ConstGeom}[1][]{
   \ifthenelse{\equal{#1}{}}       
   {K_{\First}\,}
   {K_{\First,#1}\,}
}

\newcommand{\FirstForm}[1][]{
  \ifthenelse{\equal{#1}{}}       
  {\operatorname{I_{\point}}}             
  {\operatorname{I_{#1}}}             
}
\newcommand{\GradSymbol}{\operatorname{\mathbf{\nabla}}}

\newcommand{\Lap}{\Delta}

\newcommand{\GradSurf}{\GradSymbol_{\!\SurfDomain}\!}
\newcommand{\divS}{\operatorname{div}_{\! \SurfDomain}\!}

\newcommand{\GradP}{\GradSymbol_{\!\ProjMat}\!}
\newcommand{\DivP}{\operatorname{div}_{\!\ProjMat}\!}

\newcommand{\GradC}{\GradSymbol_{\!C}\!}
\newcommand{\DivC}{\operatorname{div}_{\!C}\!}

\newcommand{\LapSurf}{\Lap_{\SurfDomain}}

\newcommand{\GradSurfConv}[1]{\GradSymbol_{#1}}

\newcommand{\gradDef}[1]{\eth_{#1}}

\newcommand{\gradFunc}[1]{\mathbf{D}_{#1}}

\newcommand{\Der}[1][]
{
  \ifthenelse{\equal{#1}{}}
  {\partial}
  {\partial_{\scriptscriptstyle{#1}}}
}
\newcommand{\DerT}{\Der_{\Time}\,}

\newcommand{\DerTot}[2][t]
{
  \ifthenelse{\equal{#2}{}}
  {\frac{d #2}{dt}}
  {\frac{d #2}{d #1}}
}
\newcommand{\Derd}[2][]{
  \ifthenelse{\equal{#1}{}}
  {\Diffsymbol{#2}}
  {\Diffsymbol_{{#1}}{#2}}
}
\newcommand{\Diffsymbol}{\operatorname{d}\!}
\newcommand{\Diff}[2][]{
  \ifthenelse{\equal{#1}{}}
  {\Diffsymbol{#2}}
  {\Diffsymbol{#2}_{{#1}}}
}
\newcommand{\DirDerSymb}{D}
\newcommand{\DirDer}[2][]
{
  \ifthenelse{\equal{#1}{}}
  {\DirDerSymb^{#2}}
  {\DirDerSymb^{#2}_{#1}}
}

\newcommand{\tr}{\operatorname{tr}}

\newcommand{\REALsymbol}{\mathbb{R}}
\newcommand{\REAL}[1][]{
  \ifthenelse{\equal{#1}{}}
  {\REALsymbol}
  {{\REALsymbol}^{#1}}
}
\newcommand{\EUCLsymbol}{\mathbb E}
\newcommand{\EUCL}[1][]{
  \ifthenelse{\equal{#1}{}}
  {\EUCLsymbol}
  {{\EUCLsymbol}^{#1}}
}
\newcommand{\NATURALsymbol}{\mathbb N}
\newcommand{\NATURAL}[1][]{
  \ifthenelse{\equal{#1}{}}
  {\NATURALsymbol}
  {{\NATURALsymbol}^{#1}}
}

\newcommand{\ABS}[2][]
{
  \ifthenelse{\equal{#1}{}}
  {\left| #2 \right|}
  {\left| #2 \right|_{#1}}
}

\newcommand{\NORM}[2][]
{
  \ifthenelse{\equal{#1}{}}
  {\left\| #2 \right\|}
  {\left\| #2 \right\|_{#1}}
}
\newcommand{\vertiii}[1]{{\left\vert\kern-0.25ex\left\vert\kern-0.25ex\left\vert #1 \right\vert\kern-0.25ex\right\vert\kern-0.25ex\right\vert}}
\newcommand{\BrNORM}[2][]
{
  \ifthenelse{\equal{#1}{}}
  {\vertiii{#2}}
  {\vertiii{#2}_{#1}}
}

\newcommand{\SCAL}[3][]
{
  \ifthenelse{\equal{#1}{}}
  {\left\langle{#2},{#3}\right\rangle}
  {\left\langle{#2},{#3}\right\rangle_{#1}}
}
\newcommand{\SCALF}[3][]
{
  \ifthenelse{\equal{#1}{}}
  {\left({#2},{#3}\right)}
  {\left({#2},{#3}\right)_{#1}}
}

\newcommand{\scalprodSurf}[3][]
{
  \ifthenelse{\equal{#1}{}}
  {\left\langle {#2},{#3} \right\rangle_{\scriptscriptstyle\SurfDomain}}
  {\left\langle {#2},{#3} \right\rangle_{{#1}}}
}

\newcommand{\InnerApprox}[2]{\left({#1}\,,\,{#2}\right)_{\meshparam}}

\newcommand{\point}[1][]
{
  \ifthenelse{\equal{#1}{}}
  {\mathbf{p}}
  {\mathbf{p}_{#1}}
}

\newcommand{\midPoint}{\mathbf{m}}

\newcommand{\RegionSymb}{R}
\newcommand{\Region}[1][]
{
  \ifthenelse{\equal{#1}{}}
  {\RegionSymb}
  {\RegionSymb_{#1}}
}

\newcommand{\Interval}[1][]
{
  \ifthenelse{\equal{#1}{}}
  {I}
  {I_{#1}}
}

\newcommand{\SurfDomainsymb}{\mathcal{S}}
\newcommand{\SurfDomain}[1][]{
  \ifthenelse{\equal{#1}{}}
  {\SurfDomainsymb}
  {{\SurfDomainsymb_{#1}}}
}
\newcommand{\tSurfDomain}[1][]{
  \ifthenelse{\equal{#1}{}}
  {\tilde{\SurfDomainsymb}}
  {\tilde{\SurfDomainsymb}_{#1}}
}
\newcommand{\SurfDomainBndsymb}{\partial\Gamma}
\newcommand{\SurfDomainBnd}[1][]{
  \ifthenelse{\equal{#1}{}}
  {\SurfDomainBndsymb}
  {\SurfDomainBndsymb_{#1}}
}
\newcommand{\ClosedSurfDomain}[1][]{
  \ifthenelse{\equal{#1}{}}
  {\Closedsymb{\SurfDomain}}
  {\Closedsymb{\SurfDomain[#1]}}
}

\newcommand{\curve}{\sigma}
\newcommand{\tcurve}{\tilde{\curve}}

\newcommand{\heightsymb}{\mathcal{H}}
\newcommand{\height}[1][]{
  \ifthenelse{\equal{#1}{}}
  {\heightsymb}
  {\heightsymb_{#1}}
}

\newcommand{\BSMsymbol}{\mathcal{B}}

\newcommand{\BSM}[1][]
{
  \ifthenelse{\equal{#1}{}}
  {\BSMsymbol}
  {\BSMsymbol_{#1}}
}

\newcommand{\Surf}{\mathcal{S}}

\newcommand{\SurfBndsymb}{\partial\Surf}
\newcommand{\SurfBnd}[1][]{
  \ifthenelse{\equal{#1}{}}
  {\SurfBndsymb}
  {\SurfBndsymb_{#1}}
}
\newcommand{\Closedsymb}[1]{\bar{#1}}
\newcommand{\ClosedSurf}[1][]{
  \ifthenelse{\equal{#1}{}}
  {\Closedsymb{\Surf}}
  {\Closedsymb{\Surf[#1]}}
}

\newcommand{\Vector}[1]{\mathbf{#1}}

\newcommand{\press}{p}

\newcommand{\densitySymb}{\rho}
\newcommand{\density}[1][]{
  \ifthenelse{\equal{#1}{}}
  {\densitySymb}
  {\densitySymb_{\scriptscriptstyle{#1}}}
}

\newcommand{\eb}{\boldsymbol{e}} 

\newcommand{\concSymbol}{\phi}
\newcommand{\CHsol}{\concSymbol}
\newcommand{\CHsolApprox}{\CHsol_{\meshparam}}
\newcommand{\pressApprox}{\press_{\meshparam}}

\newcommand{\dWellSymbol}{W}
\newcommand{\dWell}{\dWellSymbol}
\newcommand{\mobilitySymbol}{m}
\newcommand{\mobility}{\mobilitySymbol}

\newcommand{\chempot}{\mu}
\newcommand{\chempotApprox}{\mu_{\meshparam}}

\newcommand{\bendingSymb}{\kappa}
\newcommand{\bendStiff}[1][]{
  \ifthenelse{\equal{#1}{}}
  {\bendingSymb}
  {\bendingSymb_{#1}}
}

\newcommand{\bendStiffGauss}[1][]{
  \ifthenelse{\equal{#1}{}}
  {\overline{\bendingSymb}}
  {\overline{\bendingSymb}_{#1}}
}

\newcommand{\Reynolds}{\operatorname{Re}}

\newcommand{\AspRatio}{\epsilon}

\newcommand{\interfaceparam}{\AspRatio}

\newcommand{\param}{\boldsymbol{X}}

\newcommand{\update}{\boldsymbol{Y}}
\newcommand{\updateApprox}{\update_{\meshparam}}

\newcommand{\MapUsymb}{\param}
\newcommand{\MapU}[1][]
{
  \ifthenelse{\equal{#1}{}}
    {\MapUsymb}
    {\MapUsymb_{#1}}
}
\newcommand{\MapLinsymb}{F}
\newcommand{\MapLin}[1][]
{
  \ifthenelse{\equal{#1}{}}
  {\ensuremath{\MapLinsymb}}
  {\ensuremath{\MapLinsymb_{#1}}}
}

\newcommand{\MapVsymb}{\psi}
\newcommand{\MapV}[1][]
{
  \ifthenelse{\equal{#1}{}}
    {\MapVsymb}
    {\MapVsymb_{#1}}
}
\newcommand{\Transsymb}{\Phi}
\newcommand{\Trans}[1][]
{
  \ifthenelse{\equal{#1}{}}
    {\Transsymb} 
    {\Transsymb_{\scriptscriptstyle{#1}}}
}
\newcommand{\InvMapsymb}{\Psi}
\newcommand{\InvMap}[1][]
{
  \ifthenelse{\equal{#1}{}}
    {\InvMapsymb}
    {\InvMapsymb_{\scriptscriptstyle{#1}}}
}

\newcommand{\fsymb}{f}
\newcommand{\scalFun}[1][]
{
  \ifthenelse{\equal{#1}{}}
  {\fsymb}
  {\fsymb_{#1}}
}
\newcommand{\tscalFun}[1][]
{
  \ifthenelse{\equal{#1}{}}
  {\tilde{\fsymb}}
  {\tilde{\fsymb}_{#1}}
}
\newcommand{\bscalFun}[1][]
{
  \ifthenelse{\equal{#1}{}}
  {\bar{\fsymb}}
  {\bar{\fsymb}_{#1}}
}
\newcommand{\gsymb}{g}
\newcommand{\scalFung}[1][]
{
  \ifthenelse{\equal{#1}{}}
  {\gsymb}
  {\gsymb_{#1}}
}
\newcommand{\Fsymb}{F}
\newcommand{\FvecFun}[1][]
{
  \ifthenelse{\equal{#1}{}}
  {\Fsymb}
  {\Fsymb_{#1}}
}
\newcommand{\tFvecFun}[1][]
{
  \ifthenelse{\equal{#1}{}}
  {\tilde{\Fsymb}}
  {\tilde{\Fsymb}_{#1}}
}
\newcommand{\Ffunc}[2][]
{
  \ifthenelse{\equal{#1}{}}
  {\ensuremath{\Fsymb_{#2}}}
  {\ensuremath{\Fsymb^{#1}_{#2}}}
}
\newcommand{\Gsymb}{g}
\newcommand{\Gfun}[1][]
{
  \ifthenelse{\equal{#1}{}}
  {\ensuremath{\Gsymb}}
  {\ensuremath{\Gsymb_{{#1}}}}
}
\newcommand{\bGfun}[1][]
{
  \ifthenelse{\equal{#1}{}}
  {\ensuremath{\bar{\Gsymb}}}
  {\ensuremath{\bar{\Gsymb}_{{#1}}}}
}

\newcommand{\PrincipalK}[1][]
{
  \ifthenelse{\equal{#1}{}}
  {k}
  {k_{#1}}
}

\newcommand{\vecsymb}{u}
\newcommand{\vecFun}[1][]{
  \ifthenelse{\equal{#1}{}}
  {\Vector{\vecsymb}}
  {\vecsymb^{#1}}
}
\newcommand{\tvecFun}[1][]{
  \ifthenelse{\equal{#1}{}}
  {\tilde{\vecsymb}}
  {\tilde{\vecsymb}_{#1}}
}

\newcommand{\wwsymb}{w}
\newcommand{\ww}[1][]
{
  \ifthenelse{\equal{#1}{}}
  {\mathbf{\wwsymb}}
  {\wwsymb^{#1}}
}
\newcommand{\uusymb}{u}
\newcommand{\uu}[1][]
{
  \ifthenelse{\equal{#1}{}}
  {\mathbf{\uusymb}}
  {\uusymb^{#1}}
}

\newcommand{\VecFieldSymbol}{X}
\newcommand{\VecField}[1][]
{
  \ifthenelse{\equal{#1}{}}
  {\VecFieldSymbol}
  {\VecFieldSymbol^{#1}}
}
\newcommand{\VecFieldSymbolC}{Y}
\newcommand{\VecFieldYSymbol}{\boldsymbol{\VecFieldSymbolC}}
\newcommand{\VecFieldY}[1][]
{
  \ifthenelse{\equal{#1}{}}
  {\VecFieldYSymbol}
  {\VecFieldYSymbol^{#1}}
}

\newcommand{\xvsymb}{x}
\newcommand{\xv}[1][]
{
  \ifthenelse{\equal{#1}{}}
  {\mathbf{\xvsymb}}
  {\mathbf{\xvsymb}_{\scriptscriptstyle{#1}}}
}
\newcommand{\xvcomp}[1][]{
  \ifthenelse{\equal{#1}{}}
  {\xvsymb}
  {\xvsymb^{\scriptscriptstyle{#1}}}
}
\newcommand{\xcg}[1][]{
  \ifthenelse{\equal{#1}{}}
  {\xvcomp[1]}
  {\xvcomp[1]_{\scriptscriptstyle{#1}}}
}
\newcommand{\ycg}[1][]{
  \ifthenelse{\equal{#1}{}}
  {\xvcomp[2]}
  {\xvcomp[2]_{\scriptscriptstyle{#1}}}
}
\newcommand{\zcg}[1][]{
  \ifthenelse{\equal{#1}{}}
  {\xvcomp[3]}
  {\xvcomp[3]_{\scriptscriptstyle{#1}}}
}

\newcommand{\svsymb}{s}
\newcommand{\sv}[1][]{
  \ifthenelse{\equal{#1}{}}
  {\mathbf{\svsymb}}
  {\mathbf{\svsymb}_{\scriptscriptstyle{#1}}}
}
\newcommand{\svcomp}[1][]
{
  \ifthenelse{\equal{#1}{}}
  {\svsymb}
  {\svsymb^{\scriptscriptstyle{#1}}}
}
\newcommand{\xcl}[1][]{
  \ifthenelse{\equal{#1}{}}
   {\svcomp[1]}
   {\svcomp[1]_{\scriptscriptstyle{#1}}}
}
\newcommand{\ycl}[1][]{
  \ifthenelse{\equal{#1}{}}
   {\svcomp[2]}
   {\svcomp[2]_{\scriptscriptstyle{#1}}}
}
\newcommand{\zcl}[1][]{
  \ifthenelse{\equal{#1}{}}
   {\svcomp[3]}
   {\svcomp[3]_{\scriptscriptstyle{#1}}}
}

\newcommand{\ProjSymb}{\operatorname{\pi}}
\newcommand{\ProjFun}[2][]{
  \ifthenelse{\equal{#1}{}}
  {\ProjSymb\left(#2\right)}
  {\ProjSymb_{\scriptscriptstyle{#1}}\left(#2\right)}
}
\newcommand{\Prm}[1][]
{
  \ifthenelse{\equal{#1}{}}
  {\operatorname{pr}}
  {\operatorname{pr}_{\scriptscriptstyle{#1}}}
}
\newcommand{\TanPlane}[2][]
{
  \ifthenelse{\equal{#1}{}}
  {T_{\scriptscriptstyle{\point}}#2}
  {T_{\scriptscriptstyle{#1}}#2}
}

\newcommand{\SubsetSymbol}{\mathcal{U}}

\newcommand{\SubsetU}[1][]
{
  \ifthenelse{\equal{#1}{}}
  {{U}}
  {{U}_{#1}}
}
\newcommand{\SubsetV}[1][]
{
  \ifthenelse{\equal{#1}{}}
  {{V}}
  {{V}_{#1}}
}
\newcommand{\SubsetW}[1][]
{
  \ifthenelse{\equal{#1}{}}
  {{W}}
  {{W}_{#1}}
}
\newcommand{\NeighSymbol}{\mathcal{N}}
\newcommand{\Neigh}[1][]
{
  \ifthenelse{\equal{#1}{}}
  {\NeighSymbol_{\point}}
  {\NeighSymbol_{#1}}
}
\newcommand{\NeighSurf}[1][]
{
  \ifthenelse{\equal{#1}{}}
  {\SubsetSymbol_{\point}}
  {\SubsetSymbol_{#1}}
}
\newcommand{\NormSymb}{\nu}
\newcommand{\normalvec}[1][]
{
  \ifthenelse{\equal{#1}{}}
  {\boldsymbol{\NormSymb}}
  {\NormSymb_{#1}}
}\newcommand{\normalvecApprox}{\NormSymb_{\meshparam}}
\newcommand{\normalSurf}[1][]
{
  \ifthenelse{\equal{#1}{}}
  {\NormSymb}
  {\NormSymb(#1)}
}
\newcommand{\normalInterp}[1][]
{
  \ifthenelse{\equal{#1}{}}
   {\tilde{\NormSymb}}
   {\tilde{\NormSymb}_{\scriptscriptstyle{#1}}}
}

\newcommand{\normalEdge}{\mathbf{\nu}} 
\newcommand{\basisCC}{t}
\newcommand{\basisGC}{e}
\newcommand{\vecBaseGC}[1][]
{
  \ifthenelse{\equal{#1}{}}
  {\mathbf{\basisGC}}
  {\mathbf{\basisGC}_{#1}}
}
\newcommand{\vecBasePhys}[1][]
{
  \ifthenelse{\equal{#1}{}}
  {\mathbf{\basisGC}}
  {\mathbf{\basisGC}_{#1}}
}

\newcommand{\vecBaseCCcv}[1][]
{
  \ifthenelse{\equal{#1}{}}
  {\mathbf{\basisCC}}
  {\mathbf{\basisCC}_{#1}}
}

\newcommand{\tvecBaseCCcv}[1][]
{
  \ifthenelse{\equal{#1}{}}
  {\tilde{\mathbf{\basisCC}}}
  {\tilde{\mathbf{\basisCC}}_{#1}}
}
\newcommand{\hvecBaseCCcv}[1][]
{
  \ifthenelse{\equal{#1}{}}
  {\hat{\mathbf{\basisCC}}}
  {\hat{\mathbf{\basisCC}}_{#1}}
}
\newcommand{\vecBaseCCctrv}[1][]
{
  \ifthenelse{\equal{#1}{}}
  {\mathbf{\basisCC}}
  {\mathbf{\basisCC}^{#1}}
}

\newcommand{\first}[1]{
  \IfEqCase{#1}{
    {1}{\operatorname{E}}
    {2}{\operatorname{F}}
    {3}{\operatorname{G}}
  }
  [\PackageError{first}{Undefined option to first: #1}{}]%
}
\newcommand{\SecondFormSymbol}{\ensuremath{\operatorname{II}}}
\newcommand{\SecondForm}[1][]
{
  \ifthenelse{\equal{#1}{}}
  {\SecondFormSymbol_{\point}}
  {\SecondFormSymbol_{#1}}
}
\newcommand{\second}[1]{
  \IfEqCase{#1}{
    {1}{\operatorname{e}}
    {2}{\operatorname{f}}
    {3}{\operatorname{g}}
  }
  [\PackageError{first}{Undefined option to first: #1}{}]
}
\newcommand{\WeigSymbol}{\mathcal{W}}
\newcommand{\Weig}[1][]
{
  \ifthenelse{\equal{#1}{}}
  {\WeigSymbol}
  {\WeigSymbol_{#1}}
}

\newcommand{\velSymbol}{\boldsymbol{u}}
\newcommand{\vectvel}[1][]
{
   \ifthenelse{\equal{#1}{}}
   {\mathbf{\velSymbol}}
   {\mathbf{\velSymbol}(#1)}
}
\newcommand{\velVSymbol}{\boldsymbol{w}}
\newcommand{\vectvelV}[1][]
{
   \ifthenelse{\equal{#1}{}}
   {\mathbf{\velVSymbol}}
   {\mathbf{\velVSymbol}(#1)}
}
\newcommand{\vectvelApprox}[1][]
{
   \ifthenelse{\equal{#1}{}}
   {\mathbf{\velSymbol}_{\meshparam}}
   {\mathbf{\velSymbol}_{\meshparam}(#1)}
}

\newcommand{\velcompContr}[2][i]
{
   \ifthenelse{\equal{#2}{}}
   {\velSymbol^{#1}}
   {\velSymbol^{#1}(#2)}}
\newcommand{\velcompPhys}[2][i]
{
   \ifthenelse{\equal{#2}{}}
   {\velSymbol_{(#1)}}
   {\velSymbol_{(#1)}(#2)}}

\newcommand{\velSymbolRP}{v}
\newcommand{\velcompContrRP}[2][i]
{
   \ifthenelse{\equal{#2}{}}
   {\velSymbolRP^{#1}}
   {\velSymbolRP^{#1}(#2)}}
\newcommand{\velRP}[1][]
{
  \ifthenelse{\equal{#1}{}}
  {\velSymbolRP}
  {\velSymbolRP_{#1}}
}

\newcommand{\velcompApprox}[2][i]
{
   \ifthenelse{\equal{#2}{}}
   {\velSymbol^{#1)}}
   {\velSymbol^{#1}_{(#2)}}
}

\newcommand{\coord}{\mathbf{x}}

\newcommand{\VelSymbol}{U}
\newcommand{\vectVel}[1][]
{
   \ifthenelse{\equal{#1}{}}
   {\vec{\VelSymbol}}
   {\vec{\VelSymbol}(#1)}
}
\newcommand{\Velcomp}[2][i]
{
   \ifthenelse{\equal{#2}{}}
   {\VelSymbol^{#1}}
   {\VelSymbol^{#1}(#2)}
}
\newcommand{\VprimoSymbol}{\tilde{u}}
\newcommand{\Vprimo}[1][]
{
   \ifthenelse{\equal{#1}{}}
   {\VprimoSymbol}
   {\VprimoSymbol(#1)}
}
\newcommand{\VprimoComp}[2][i]
{
   \ifthenelse{\equal{#2}{}}
   {\VprimoSymbol^{#1}}
   {\VprimoSymbol^{#1}(#2)}
}
\newcommand{\ttvelSymbol}{\tilde{\Mymathbb{u}}}
\newcommand{\ttvel}[1][]
{
   \ifthenelse{\equal{#1}{}}
   {\mathbf{\ttvelSymbol}}
   {\mathbf{\ttvelSymbol}(#1)}
}
\newcommand{\ttvelComp}[2][i]
{
   \ifthenelse{\equal{#2}{}}
   {\ttvelSymbol^{#1}}
   {\ttvelSymbol^{#1}(#2)}
}

\newcommand{\MatAlphaSymbol}{\mathbb{A}}
\newcommand{\MatAlpha}[1][]{%
  \ifthenelse{\equal{#1}{}}
  {\MatAlphaSymbol}
  {\MatAlphaSymbol_{#1}}
}

\newcommand{\QSymbol}{q}
\newcommand{\Qdisch}[1][]
{
   \ifthenelse{\equal{#1}{}}
   {\mathbf{\QSymbol}}
   {\mathbf{\QSymbol}(#1)}
}
\newcommand{\Qcomp}[2][i]
{
   \ifthenelse{\equal{#2}{}}
   {\QSymbol^{#1}}
   {\QSymbol^{#1}(#2)}
}
\newcommand{\Qvect}[1][]
{
   \ifthenelse{\equal{#1}{}}
   {\mathbf{\QSymbol}}
   {\mathbf{\QSymbol}(#1)}
}
\newcommand{\FricSymbol}{f}
\newcommand{\vectFric}[1][]
{
   \ifthenelse{\equal{#1}{}}
   {\mathbf{\FricSymbol}}
   {\mathbf{\FricSymbol}_{\scriptscriptstyle{#1}}}
}
\newcommand{\Friccomp}[2][i]
{
   \ifthenelse{\equal{#2}{}}
   {\FricSymbol_{#1}}
   {\FricSymbol_{#1}(#2)}}
\newcommand{\BFsymbol}{\tau}
\newcommand{\BottomFriction}[1][]{
  \ifthenelse{\equal{#1}{}}
  {\BFsymbol_{b}}
  {\BFsymbol_{b}^{#1}}
}

\newcommand{\ProjMatSymb}{\boldsymbol{P}}
\newcommand{\ProjMat}[1][]{
  \ifthenelse{\equal{#1}{}}
  {\ProjMatSymb}
  {\ProjMatSymb_{#1}}
}
\newcommand{\IDSymbol}{\boldsymbol{I}}
\newcommand{\IDtens}[1][]{
  \ifthenelse{\equal{#1}{}}
  {\IDSymbol}
  {\IDSymbol(#1)}
}
\newcommand{\IDMat}{\IDSymbol}

\newcommand{\tensSymbol}{\boldsymbol{T}}
\newcommand{\tenscompSymbol}{\tau}
\newcommand{\tens}[1][]{
  \ifthenelse{\equal{#1}{}}
  {\tensSymbol}
  {\tensSymbol(#1)}
}
\newcommand{\tenscomp}[2][ij]
{
  \ifthenelse{\equal{#2}{}}
  {\tenscompSymbol^{#1}}
  {\tenscompSymbol^{#1}(#2)}
}
\newcommand{\tensrow}[2][i]
{
  \ifthenelse{\equal{#2}{}}
  {\tensSymbol^{(#1)}}
  {\tensSymbol^{(#1)}(#2)}
}
\newcommand{\TensSymbol}{\mathbf{T}}
\newcommand{\Tens}[1][]{
  \ifthenelse{\equal{#1}{}}
  {\TensSymbol}
  {\TensSymbol_{#1}}
}

\newcommand{\TensCompSymbol}{\TensSymbol}
\newcommand{\TensComp}[2][ij]
{
  \ifthenelse{\equal{#2}{}}
  {\TensCompSymbol^{#1}}
  {\TensCompSymbol^{#1}(#2)}
}

\newcommand{\tensPrimoSymbol}{\tilde{\mathbf{\tau}}}
\newcommand{\tensPrimo}[1][]{
  \ifthenelse{\equal{#1}{}}
  {\tensPrimoSymbol}
  {\tensPrimoSymbol(#1)}
}
\newcommand{\tensPrimoCompSymbol}{\tensPrimoSymbol}
\newcommand{\tensPrimoComp}[2][ij]
{
  \ifthenelse{\equal{#2}{}}
  {\tensPrimoCompSymbol^{#1}}
  {\tensPrimoCompSymbol^{#1}(#2)}
}
\newcommand{\MCxl}[1][]{\ifthenelse{\equal{#1}{}}{h_{(1)}}{h_{(1),#1}}} 
\newcommand{\MCyl}[1][]{\ifthenelse{\equal{#1}{}}{h_{(2)}}{h_{(2),#1}}} 
\newcommand{\MCzl}[1][]{\ifthenelse{\equal{#1}{}}{h_{(3)}}{h_{(3),#1}}}

\newcommand{\MPsymb}{h}

\newcommand{\smallestH}[1][]
{
  \ifthenelse{\equal{#1}{}}
    {{l}}
    {{l}_{#1}}
}
\newcommand{\meshparam}[1][]
{
  \ifthenelse{\equal{#1}{}}
    {\MPsymb}
    {\MPsymb_{\scriptscriptstyle{#1}}}
}
\newcommand{\InradiusSymbol}{r}
\newcommand{\Inradius}[1][]
{
  \ifthenelse{\equal{#1}{}}
  {\InradiusSymbol}
  {\InradiusSymbol_{\scriptscriptstyle{#1}}}
}
\newcommand{\Tsymb}{\mathcal{T}}
\newcommand{\Triang}[1][]
{
  \ifthenelse{\equal{#1}{}}
    {\Tsymb}
    {\Tsymb_{#1}}
}
\newcommand{\TriangH}[1][]
{
  \ifthenelse{\equal{#1}{}}
    {\Tsymb_{\meshparam}}
    {\Tsymb_{#1}}
}
\renewcommand{\Triang}{\TriangH}
\newcommand{\Edgesymb}{\sigma}
\newcommand{\Edge}[1][]{
  \ifthenelse{\equal{#1}{}}
    {\Edgesymb}
    {\Edgesymb_{#1}}
}
\newcommand{\EdgeH}[1][]{
  \ifthenelse{\equal{#1}{}}
    {\Edgesymb_{\meshparam}}
    {\Edgesymb_{\meshparam,#1}}
}
\newcommand{\NEdge}[1][]{
  \ifthenelse{\equal{#1}{}}
    {N_{\Edgesymb}}
    {N_{\Edgesymb({#1})}}
}
\newcommand{\Cellsymb}{T}
\newcommand{\Cell}[1][]{
  \ifthenelse{\equal{#1}{}}
    {\Cellsymb}
    {\Cellsymb_{#1}}
}
\newcommand{\tCell}[1][]{
  \ifthenelse{\equal{#1}{}}
    {\tilde{\Cellsymb}}
    {\tilde{\Cellsymb}_{#1}}
}
\newcommand{\CellH}[1][]{
  \ifthenelse{\equal{#1}{}}
    {\Cellsymb_{\meshparam}}
    {\Cellsymb_{\meshparam,#1}}
}
\newcommand{\areaSymb}{\mathcal{A}}
\newcommand{\CellArea}[1][]
{
  \ifthenelse{\equal{#1}{}}
    {\areaSymb_{\Cell}}
    {\areaSymb_{#1}}
}
\newcommand{\CellHArea}[1][]
{
  \ifthenelse{\equal{#1}{}}
    {\areaSymb_{\CellH}}
    {\areaSymb_{\meshparam,#1}}
}

\newcommand{\NCell}[1][]{
  \ifthenelse{\equal{#1}{}}
    {N_{\Cellsymb}}
    {N_{\Cellsymb({#1})}}
}

\newcommand{\lengthSymb}{l}

\newcommand{\edgeLength}[1][]
{
  \ifthenelse{\equal{#1}{}}
  {\lengthSymb_{\Edge}}
  {\lengthSymb_{#1}}
}
 
\newcommand{\edgeHLength}[1][]
{
  \ifthenelse{\equal{#1}{}}
  {\lengthSymb_{\EdgeH}}
  {\lengthSymb_{\meshparam,#1}}
}

\newcommand{\Sourcesymb}{\mathbf{S}}
\newcommand{\Source}[1][]
{
  \ifthenelse{\equal{#1}{}}
    {\Sourcesymb}
    {\Sourcesymb_{#1}}
}

\newcommand{\SourceEdge}[1][]{
  \ifthenelse{\equal{#1}{}}
    {\Sourcesymb_{ij}}
    {\Sourcesymb_{#1}}
}

\newcommand{\FluxEdgesymb}{\mathbf{F}}
\newcommand{\FluxEdge}[1][]{
  \ifthenelse{\equal{#1}{}}
    {\FluxEdgesymb_{ij}}
    {\FluxEdgesymb_{#1}}
}

\newcommand{\numFlux}[1][]
{
  \ifthenelse{\equal{#1}{}}
  {\tilde{\fsymb}}
  {\tilde{\fsymb}_{#1}}
}

\newcommand{\FluxFuncNormSymbol}{\mathbf{F}}
\newcommand{\FluxFuncNorm}[1][]{
  \ifthenelse{\equal{#1}{}}
    {\FluxFuncNormSymbol^{\normalEdge}}
    {\FluxFuncNormSymbol^{\normalEdge}_{#1}}    
}

\newcommand{\JacobianSymbol}{\mathbf{A}}
\newcommand{\Jacobian}[1][]{
  \ifthenelse{\equal{#1}{}}
    {\JacobianSymbol}
    {\JacobianSymbol_{#1}}    
}
\newcommand{\EValSymbol}{\lambda}
\newcommand{\EVal}[1][]{
  \ifthenelse{\equal{#1}{}}
  {\EValSymbol}
  {\EValSymbol_{#1}}
}
\newcommand{\EVecSymbol}{\mathbf{r}}
\newcommand{\EVec}[2][]{
  \ifthenelse{\equal{#1}{}}
  {\EVecSymbol^{(#2)}}
  {\EVecSymbol^{(#2)}_{#1}}
}

\newcommand{\midPointEdge}[1][]
{
  \ifthenelse{\equal{#1}{}}
  {\midPoint_{\scriptscriptstyle\Edge}}
  {\midPoint_{\scriptscriptstyle\Edge[#1]}}
}
\newcommand{\gpPointEdgeDG}[1][]
{
  \ifthenelse{\equal{#1}{}}
  {\point_{\scriptscriptstyle\Edge}}
  {\point_{\scriptscriptstyle\Edge,#1}}
}
\newcommand{\midPointCell}[1][]
{
  \ifthenelse{\equal{#1}{}}
  {\midPoint_{\scriptscriptstyle\Cell}}
  {\midPoint_{\scriptscriptstyle\Cell[#1]}}
}

\newcommand{\energySymbol}{\mathcal{F}}
\newcommand{\energy}[1][]
{
  \ifthenelse{\equal{#1}{}}
  {\energySymbol}
  {\energySymbol_{#1}}
}

\newcommand{\energyHelfr}{\energy[H]}
\newcommand{\energyGL}{{\energy[GL]}}
\newcommand{\energyKin}{{\energy[K]}}
\newcommand{\Lagrangian}{\mathcal{L}}
\newcommand{\DissPotential}{\mathcal{D}}

\newcommand{\speedRS}[2][]
{
  \ifthenelse{\equal{#2}{}}
  {S_{#2}}
  {S_{#2}^{#1}}
}
\newcommand{\Interpolant}{I_{\meshparam}}

\newcommand{\ContSymbol}{C}
\newcommand{\Cont}[1][]{
  \ifthenelse{\equal{#1}{}}
  {\ContSymbol^{0}}
  {\ContSymbol^{#1}}
}
\newcommand{\Cinf}[1][]{
  \ifthenelse{\equal{#1}{}}
  {\ContSymbol^{\infty}}
  {\ContSymbol^{\infty}(#1)}
}
\newcommand{\SobSymbol}{W}
\newcommand{\HilbSymbol}{H}
\newcommand{\Hilb}[1][]{
  \ifthenelse{\equal{#1}{}}
  {\HilbSymbol^{1}}
  {\SobSymbol^{1}_{#1}}  
}
\newcommand{\Sob}[2][]{
  \ifthenelse{\equal{#1}{}}
  {\HilbSymbol^{#2}}
  {\SobSymbol^{#2,#1}}  
}

\newcommand{\LspaceSymb}{L}
\newcommand{\Lspace}[1][]{
  \ifthenelse{\equal{#1}{}}
  {\LspaceSymb^{2}}
  {\LspaceSymb^{#1}}  
}

\newcommand{\TestSpSymbol}{{V}}
\newcommand{\TestSpace}[1][]{
  \ifthenelse{\equal{#1}{}}
  {\TestSpSymbol({\SurfDomain})}
  {\TestSpSymbol_{#1}({\SurfDomain})}
}

\newcommand{\TestSpaceEmbedded}[1][]{
  \ifthenelse{\equal{#1}{}}
  {\TestSpSymbol(\TriangH(\SurfDomain))}
  {\TestSpSymbol_{#1}(\TriangH(\SurfDomain))}
}
\newcommand{\TestSpaceIntrinsic}[1][]{
  \ifthenelse{\equal{#1}{}}
  {\TestSpSymbol(\Triang(\SurfDomain))}
  {\TestSpSymbol_{#1}(\Triang(\SurfDomain))}
}
\newcommand{\TestSpaceChart}[1][]{
  \ifthenelse{\equal{#1}{}}
  {\TestSpSymbol(\Triang(\SubsetU))}
  {\TestSpSymbol_{#1}(\Triang(\SubsetU))}
}
\newcommand{\TestSpaceCell}[1][]{
  \ifthenelse{\equal{#1}{}}
  {\TestSpSymbol_{\meshparam}(\Cell)}
  {\TestSpSymbol_{\meshparam}(\Cell_{#1})}
}

\newcommand{\TestSpaceApprox}[1][]{
  \ifthenelse{\equal{#1}{}}
  {\TestSpSymbol_{\meshparam}}
  {\TestSpSymbol_{\meshparam}(#1)}
}
\newcommand{\TestSpaceVec}[1][]{
  \ifthenelse{\equal{#1}{}}
  {\mathbf{\TestSpSymbol}_{\meshparam}}
  {\mathbf{\TestSpSymbol}_{\meshparam}(#1)}
}

\newcommand{\TestSpGammaSymbol}{{P}}
\newcommand{\TestSpaceGamma}[1][]{
  \ifthenelse{\equal{#1}{}}
  {\TestSpGammaSymbol_{\meshparam}}
  {\TestSpGammaSymbol_{\meshparam}(#1)}
}

\newcommand{\TestSpCHSymbol}{{W}}
\newcommand{\TestSpaceCH}[1][]{
  \ifthenelse{\equal{#1}{}}
  {\TestSpCHSymbol_{\meshparam}}
  {\TestSpCHSymbol_{\meshparam}(#1)}
}

\newcommand{\TestSpPressSymbol}{Q}
\newcommand{\TestSpacePress}[1][]{
  \ifthenelse{\equal{#1}{}}
  {\TestSpPressSymbol_{\meshparam}}
  {\TestSpPressSymbol_{\meshparam}(#1)}
}

\newcommand{\GaussCurv}{\mathcal{K}}
\newcommand{\calH}{\mathcal{H}}
\newcommand{\meanCurv}[1][]
{
  \ifthenelse{\equal{#1}{}}
  {\calH}
  {\calH_{#1}}
}
\newcommand{\hmeanCurv}[1][]
{
  \ifthenelse{\equal{#1}{}}
  {\hat{\calH}}
  {\hat{\calH}_{#1}}
}
\newcommand{\meanCurvApprox}{\meanCurv[\meshparam]}
\newcommand{\shapeOp}{\mathcal{B}}
\newcommand{\AngleSymbol}{\theta}
\newcommand{\DevAngle}[1][]
{
  \ifthenelse{\equal{#1}{}}
  {\AngleSymbol}
  {\AngleSymbol_{\scriptscriptstyle{#1}}}
}
\newcommand{\relheightSymb}{\pi}
\newcommand{\relheight}[1][]
{
  \ifthenelse{\equal{#1}{}}
  {\relheightSymb_{\scriptscriptstyle{\SurfDomain}}}
  {\relheightSymb_{\scriptscriptstyle{#1}}}
}

\newcommand{\KinCoef}{\gamma}

\newcommand{\StressTensComp}{\sigma}
\newcommand{\StressTens}{\boldsymbol{\StressTensComp}}

\newcommand{\DiffEig}[1][]
{ \ifthenelse{\equal{#1}{}}
  {d}
  {d_{#1}}
}

\newcommand{\BilinearStiffSymbol}{a}
\newcommand{\BilinearStiff}[3][]
{
  \ifthenelse{\equal{#1}{}}
  {\BilinearStiffSymbol(#2,#3)}    
  {\BilinearStiffSymbol_{#1}(#2,#3)}    
}
\newcommand{\BilinearAdvSymbol}{b}
\newcommand{\BilinearAdv}[3][]
{
  \ifthenelse{\equal{#1}{}}
  {\BilinearAdvSymbol(#2,#3)}    
  {\BilinearAdvSymbol_{#1}(#2,#3)}    
}
\newcommand{\BilinearMassSymbol}{m}
\newcommand{\BilinearMass}[3][]
{
  \ifthenelse{\equal{#1}{}}
  {\BilinearMassSymbol(#2,#3)}    
  {\BilinearMassSymbol_{#1}(#2,#3)}    
}
\newcommand{\BilinearReactSymbol}{c}
\newcommand{\BilinearReact}[3][]
{
  \ifthenelse{\equal{#1}{}}
  {\BilinearReactSymbol(#2,#3)}    
  {\BilinearReactSymbol_{#1}(#2,#3)}    
}

\newcommand{\TestSymbol}{v}
\newcommand{\Test}[1][]
{
  \ifthenelse{\equal{#1}{}}
  {\TestSymbol}  
  {\TestSymbol_{\scriptscriptstyle{#1}}}
}

\newcommand{\nNodes}[1][]
{
  \ifthenelse{\equal{#1}{}}
  {N^{\scriptscriptstyle{dof}}}
  {N^{\scriptscriptstyle{dof}}_{\scriptscriptstyle{#1}}}
}

\newcommand{\TestApprox}[1][]
{
  \ifthenelse{\equal{#1}{}}
  {\TestSymbol_{\scriptscriptstyle{\meshparam}}}  
  {\TestSymbol_{\scriptscriptstyle{\meshparam,#1}}}
}

\newcommand{\bTestApprox}[1][]
{
  \ifthenelse{\equal{#1}{}}
  {\bar{\TestSymbol}_{\scriptscriptstyle{\meshparam}}}  
  {\bar{\TestSymbol}_{\scriptscriptstyle{\meshparam,#1}}}
}

\newcommand{\TestZSymbol}{\mathbf{Z}}
\newcommand{\TestZ}[1][]
{
  \ifthenelse{\equal{#1}{}}
  {\TestZSymbol}  
  {\TestZSymbol_{\scriptscriptstyle{#1}}}
}

\newcommand{\TestZApprox}[1][]
{
  \ifthenelse{\equal{#1}{}}
  {\TestZSymbol_{\scriptscriptstyle{\meshparam}}}  
  {\TestZSymbol_{\scriptscriptstyle{\meshparam,#1}}}
}

\newcommand{\TestHSymbol}{{h}}
\newcommand{\TestH}[1][]
{
  \ifthenelse{\equal{#1}{}}
  {\TestHSymbol}  
  {\TestHSymbol_{\scriptscriptstyle{#1}}}
}

\newcommand{\TestHApprox}[1][]
{
  \ifthenelse{\equal{#1}{}}
  {\TestHSymbol_{\scriptscriptstyle{\meshparam}}}  
  {\TestHSymbol_{\scriptscriptstyle{\meshparam,#1}}}
}

\newcommand{\TestWSymbol}{w}
\newcommand{\TestW}[1][]
{
  \ifthenelse{\equal{#1}{}}
  {\TestWSymbol}  
  {\TestWSymbol_{\scriptscriptstyle{#1}}}
}

\newcommand{\TestWApprox}[1][]
{
  \ifthenelse{\equal{#1}{}}
  {\TestWSymbol_{\scriptscriptstyle{\meshparam}}}  
  {\TestWSymbol_{\scriptscriptstyle{\meshparam,#1}}}
}
\newcommand{\bTestWApprox}[1][]
{
  \ifthenelse{\equal{#1}{}}
  {\bar{\TestWSymbol}_{\scriptscriptstyle{\meshparam}}}  
  {\bar{\TestWSymbol}_{\scriptscriptstyle{\meshparam,#1}}}
}

\newcommand{\TestCHSymbol}{\psi}
\newcommand{\TestCHApprox}[1][]
{
  \ifthenelse{\equal{#1}{}}
  {\TestCHSymbol_{\scriptscriptstyle{\meshparam}}}  
  {\TestCHSymbol_{\scriptscriptstyle{\meshparam,#1}}}
}
\newcommand{\TestCH}[1][]
{
  \ifthenelse{\equal{#1}{}}
  {\TestCHSymbol}  
  {\TestCHSymbol_{\scriptscriptstyle{#1}}}
}
\newcommand{\TestCHmuSymbol}{\xi}
\newcommand{\TestCHmuApprox}[1][]
{
  \ifthenelse{\equal{#1}{}}
  {\TestCHmuSymbol_{\scriptscriptstyle{\meshparam}}}  
  {\TestCHmuSymbol_{\scriptscriptstyle{\meshparam,#1}}}
}

\newcommand{\TestVecSymbol}{\boldsymbol{v}}
\newcommand{\TestVecApprox}[1][]
{
  \ifthenelse{\equal{#1}{}}
  {\TestVecSymbol_{\scriptscriptstyle{\meshparam}}}  
  {\TestVecSymbol_{\scriptscriptstyle{\meshparam,#1}}}
}
\newcommand{\TestPressSymbol}{q}
\newcommand{\TestPress}[1][]
{
  \ifthenelse{\equal{#1}{}}
  {\TestPressSymbol}  
  {\TestPressSymbol_{\scriptscriptstyle{#1}}}
}
\newcommand{\TestPressApprox}[1][]
{
  \ifthenelse{\equal{#1}{}}
  {\TestPressSymbol_{\scriptscriptstyle{\meshparam}}}  
  {\TestPressSymbol_{\scriptscriptstyle{\meshparam,#1}}}
}

\newcommand{\viscositySymb}{\eta}
\newcommand{\viscosity}[1][]
{
  \ifthenelse{\equal{#1}{}}
  {\viscositySymb}
  {\viscositySymb_{\scriptscriptstyle{#1}}}
}

\newcommand{\ResidualSymbol}{R}
\newcommand{\Residual}[1][]
{
  \ifthenelse{\equal{#1}{}}
  {\ResidualSymbol}
  {\ResidualSymbol_{#1}}
}

\newcommand{\curvature}[1][]
{
  \ifthenelse{\equal{#1}{}}
  {\kappa}
  {\kappa_{#1}}
}

\newcommand{\QuadRule}[2][]
{
  \ifthenelse{\equal{#1}{}}
  {Q(#2)}
  {Q_{#1}(#2)}
}

\newcommand{\mfrak}{\mathfrak{m}}

\newcommand{\ie}{i.\,e.}

\newcommand{\formComma}{\,\text{,}}
\newcommand{\formPeriod}{\,\text{.}}

\newcommand{\xb}{\boldsymbol{x}}

\newcommand{\wb}{\vectvelV}
\newcommand{\ub}{\vectvel}
\newcommand{\vb}{\boldsymbol{\upsilon}}

\usepackage{cancel}
\usepackage{changes}
\newcommand{\soutthick}[1]{%
    \renewcommand{\ULthickness}{1.0pt}%
       \sout{#1}%
    \renewcommand{\ULthickness}{.4pt}
}
\newcommand{\stkout}[1]{
  \ifmmode\text{\soutthick{\ensuremath{#1}}}\else\soutthick{#1}\fi
}
\setdeletedmarkup{\stkout{#1}}

\title{Derivation and simulation of a two-phase fluid deformable surface model}
\author{
Elena Bachini$^1$
\footnote{Now at: Department of Mathematics ``Tullio Levi-Civita'', University of Padua, Italy} 
\and Veit Krause$^1$
\and Ingo Nitschke$^1$  
\and Axel Voigt$^{1,2,3}$\footnote{Email address for correspondence: axel.voigt@tu-dresden.de}\\[1em]
$^1$Institute of Scientific Computing, TU Dresden, Germany\\[0.5em]
$^2$Center for System Biology Dresden, Germany\\[0.5em]
$^3$Cluster of Excellence, Physics of Life, TU Dresden, Germany
}

\date{}

\begin{document}
\maketitle

\begin{abstract}
We consider two-phase fluid deformable surfaces as model systems for biomembranes. Such surfaces are modeled by incompressible surface Navier-Stokes-Cahn-Hilliard-like equations with bending forces. We derive this model using the Lagrange-D'Alembert principle considering various dissipation mechanisms. The highly nonlinear model is solved numerically to explore the tight interplay between surface evolution, surface phase composition, surface curvature and surface hydrodynamics. It is demonstrated that hydrodynamics can enhance bulging and furrow formation, which both can further develop to pinch-offs. The numerical approach builds on a Taylor-Hood element for the surface Navier-Stokes part, a semi-implicit approach for the Cahn-Hilliard part, higher order surface parametrizations,  appropriate approximations of the geometric quantities, and mesh redistribution. We demonstrate convergence properties that are known to be optimal for simplified sub-problems.
\end{abstract}

\section{Introduction}
\label{sec:intro}

Coexisting fluid domains in biomembranes and in their model systems is an intensively studied field of research. With the possibility to visualize coexisting lipid phases in model membranes by high-resolution fluorescence imaging \cite{BHW_N_2003} a strong correlation between domain composition and local membrane curvature can be established. These findings have supported earlier membrane models \cite{julicher1996shape,jiang2000phase,kumar2001budding} and initiated a wealth of theoretical, numerical, and experimental studies to understand the relation between composition and curvature \cite{VEATCH20033074,BHW_N_2003,BDWJ_BPJ_2005,veatch2007critical,WD_JMB_2008,LRV_PRE_2009,art:elliott2010surface,EFDG_ARBP_2010,GKRR_MMMAS_2016,GS_CF_2018,Zimmermannetal_CMAME_2019}. This relation has strong biological implications \cite{MG_N_2005}. We refer to \cite{VEATCH20033074,Deserno_2015,art:Lipowsky2021SM} for reviews on the subject. While most studies address multicomponent giant unilaminar vesicles (GUV) and consider phase separation and coarsening on spherical shapes, more recently multicomponent scaffolded lipid vesicles (SLV) have been considered. In these systems, non-spherical shapes are stabilized and the effect of spatially varying curvature on phase separation and coarsening has been considered \cite{fonda2018interface,art:Fondaetal2019,RFGK_NC_2020}. These results demonstrate the strong influence of curvature on the spatial arrangement of the lipid phases. On the other hand also the shape evolution is considered. Coexisting fluid domains can lead to bulging and budding events \cite{BHW_N_2003,LRV_PRE_2009,art:elliott2010surface}, indicating also the strong influence of composition on the evolving shape. Thus, the correlation between domain composition and membrane curvature acts in both directions. These effects can essentially be modeled by the J\"ulicher-Lipowsky model \cite{julicher1996shape} or appropriate phase field approximations of it \cite{LRV_PRE_2009,Hauser2013}. These models essentially extend the classical Helfrich model to multiple components. Instead of a simple $L^2$-gradient flow of the Helfrich energy, with possible constraints, the model allows for dissipation through the simultaneous evolution of the shape and the lipid phases on the surface. Such models can be embedded in bulk flows and serve as interfacial conditions \cite{SOHN2010119,ZHANG2022110815}. However, all these studies neglect the effect of surface viscosity. Surface viscosity has been shown to be a key property of biomembranes and their model systems controlling remodelling and with it, the coarsening process, see \cite{FaiZietal_BPJ_2022} for discussions and measurements. Already for flat membranes, it has been shown that experimental results on the coarsening rate of these fluid domains can only be quantitatively reproduced by simulations if the fluid properties of the membrane are taken into account \cite{FHH_JCP_2010,CB_JCP_2011}. Considering these effects on curved surfaces involves several modelling and numerical subtleties, due to the increased coupling between local curvature and surface fluid velocity. Although there are models in the literature dealing with the coupling between surface flow on curved surfaces and the surrounding bulk flow \cite{Barrettetal_PRE_2015,woodhouse_goldstein_2012,reuther2016incompressible}, they do not consider coexisting fluid domains. Furthermore, it has been shown that the effect of the bulk fluid can be neglected if the {Saffman-Delbr\"uck} number 
\cite{SD_PNAS_1975}, which defines a hydrodynamic length relating the viscosity of the membrane and the surrounding bulk fluid, is large. In this case the hydrodynamics is effectively 2D on spatial scales smaller than the Saffmann-Delbr\"uck number \cite{PhysRevE.81.011905}. We follow this simplification and consider only surface two-phase flow problems. On stationary surfaces these problems are addressed in \cite{Nitschke_2012,ABWPK_PRE_2019,Olshanskiietal_VJM_2022,art:BKV23} and essentially reveal an enhanced coarsening rate due to hydrodynamic effects, similar to the situation in flat space \cite{FHH_JCP_2010,CB_JCP_2011}. To consider evolving surfaces under the influence of surface viscosity requires models for so-called fluid deformable surfaces. They exhibit a solid-fluid duality, while they store elastic energy when stretched or bent, as solid shells, they flow as viscous two-dimensional fluids under in-plane shear. This duality has several consequences: it establishes a tight interplay between tangential flow and surface deformation. In the presence of curvature, any shape change is accompanied by a tangential flow and, vice-versa, the surface deforms due to tangential flow. Models describing this interplay between curvature and surface flow have been introduced and numerically solved in \cite{art:torres-sanchez_millan_arroyo_2019,reuther2020numerical,art:Krause22}. We here extend these models to two-phase flows by combining fluid deformable surfaces with a phase field approximation of the J\"ulicher-Lipowsky model. The model is systematically derived by a Lagrange-D'Alembert principle, accounting for dissipation by domain coarsening, shape evolution and surface viscosity. We relate the resulting model to known simplified models in the literature. We provide a detailed description of the numerical approach, which uses surface finite elements (SFEM) \cite{art:Dziuk2013,NNV_JCP_2019} and builds on previous developments for one-component fluid deformable surfaces \cite{art:Krause22} and surface Navier-Stokes-Cahn-Hilliard-like models on stationary surfaces \cite{art:BKV23}. The algorithm is used to demonstrate the strong interplay between composition, curvature and hydrodynamics and its implications for bulging and budding processes. Briefly, surface hydrodynamics enhances the onset of such shape deformations and possible resulting topological changes, which has strong biological implications, e.g., in the case of endocytosis and exocytosis \cite{Kasonenetal_NRMCB_2018,AlIzzi2020}. 

The paper is structured as follows. In Section \ref{sec:model} we introduce the used notation necessary to formulate the surface model, briefly describe the model derivation, formulate the two-phase fluid deformable surface model, and relate the equations to known simplified models. Details are considered in the Appendices. In Section \ref{sec:discrete} we discuss the considered numerical approach. In Section \ref{sec:results} we provide all used parameters, consider convergence studies, and explore the parameter space and the implications of the coupling between composition, shape and surface flow. In Section \ref{sec:conclusion} we draw conclusions.

\section{Continuous model}
\label{sec:model}

We derive the full model for two-phase fluid deformable surfaces by applying the Lagrange-D'Alembert principle. We first introduce the necessary notation, then motivate the use of the Lagrange-D'Alembert principle, introduce all ingredients, and derive the model. Finally, we relate the model to known simplifications in the literature.

\subsection{Notation}

We consider a time dependent smooth and oriented surface $\SurfDomain = \SurfDomain(t)$ without boundary.
Related to $\SurfDomain$, we denote the outward pointing surface normal $\normalvec$, the surface projection $\ProjMat=\IDMat-\normalvec\otimes\normalvec$, the shape operator $\shapeOp= -\GradP\normalvec$, and the mean curvature $\meanCurv= \operatorname{tr}\shapeOp$. Note that, under these definitions the unit sphere has negative mean curvature $\meanCurv=-2$.  
Let $\phi$ be a continuously differentiable scalar field, $\vectvel$ a continuously differentiable $\REAL^3$-vector field, and $\StressTens$ a continuously differentiable $\REAL^{3\times3}$-tensor field defined on $\SurfDomain$. We define the surface tangential gradient $\GradP$ as in~ \cite{jankuhn2018} and the componentwise surface gradient $\GradC$ as in~\cite{Nitschke2022GoI}, namely:
\begin{align*}
  \GradP\phi &= \ProjMat\nabla\phi^e\,, \\
  \GradP\vectvel &= \ProjMat\nabla\vectvel^e\ProjMat \,, \\
  \GradC\StressTens &= \nabla\StressTens^e \ProjMat\,,
\end{align*}
where $\phi^e$, $\vectvel^e$ and $\StressTens^e$ are arbitrary smooth extensions of $\phi$, $\vectvel$ and $\StressTens$ in the normal direction and $\nabla$ is the gradient of the embedding space $\REAL^3$. The fields $\GradP\phi, \GradP\vectvel$ are purely tangential vector and tensor fields, respectively.  We define the corresponding divergence operators for a vector field $\vectvel$ and a tensor field $\StressTens$ by:
\begin{align*}
  \DivP\vectvel  &= \operatorname{tr}(\GradP\vectvel)\,, \\
  \DivC(\StressTens\ProjMat) &= \operatorname{tr}\GradC(\StressTens\ProjMat)\,,
\end{align*}
where $\operatorname{tr}$ is the trace operator. The divergence of the tensor $\StressTens$,  $\DivC\StressTens$, leads to a non-tangential vector field even if $\StressTens$ is a tangential tensor field. Let $\GradSurf$ be the gradient with respect to the covariant derivative on $\SurfDomain$, as used in \cite{reuther2020numerical}. This operator is defined for scalar fields and tangential vector fields and relates to the tangential operators by $\GradP\phi=\GradSurf\phi$ and
${\DivP\vectvel = \divS(\ProjMat\vectvel)-(\vectvel\cdot\normalvec)\meanCurv}$, respectively.
We further clarify the relation between the above operators in Appendix \ref{app:symbols}.

The surface $\SurfDomain$ is given by a parametrization $\param$. The material on the surface is described by a material parametrization $\param_\mfrak$, as in~\cite{nitschke2022TimeDerivative}. Both parametrizations relate to each other by $\partial_t\param_\mfrak\cdot\normalvec = \partial_t\param \cdot \normalvec$, which leads to a Lagrangian perspective in normal direction. In tangential direction the surface and the material can move independently. We define the material velocity by $\vectvel\coloneqq \partial_t\param_\mfrak$ and the relative material velocity by $\vectvelV := \vectvel - \partial_t\param$. The relative material velocity is a pure tangential vector field. Additional information and the relation to the time derivative used can be found in Appendix \ref{sec:material_derivative}. 

\subsection{Model derivation by Lagrange-D'Alembert principle}

The Lagrange-D'Alembert principle has been explained in \cite{marsden2013introduction}. The approach is a combination of the Lagrange and the Onsager variational principles. The concept of the Lagrange principle has been introduced in \cite{marsden2001discrete,marsden2013introduction,hairer2006geometric} and it models the interplay of kinetic and potential energies under total energy conservation. On the opposite, the concept of the Onsager variational principle models dissipative systems, where potential energy decreases under a dissipation potential, e.g. $L^2$-gradient flows. Within our context, the Onsager principle is used in \cite{art:torres-sanchez_millan_arroyo_2019} to derive fluid deformable surfaces with the surface Stokes model. 

We consider a density function $\CHsol$ that describes the two-phases on the surface $\SurfDomain$,
e.g., liquid-ordered and disordered phases, where one phase is represented by $\phi=1$, the other one by $\phi=-1$, and we assume a mixture of both if $\CHsol \in (-1,1)$. We consider a conserved evolution for $\CHsol$ with respect to the following energies.

The first is a Ginzburg-Landau energy modeling phase-separation on $\SurfDomain$:
\begin{equation}
    \label{eq:GL}
    \energyGL = \int_{\SurfDomain}
    \tcurve \left( \frac{\interfaceparam}{2}\NORM{\GradSurf\CHsol}^2
    +\frac{1}{\interfaceparam}\dWell(\CHsol) \right) \Diff\SurfDomain\,,
\end{equation}
where $\tcurve>0$ is a rescaled line tension with $\tcurve = 3 / (2 \sqrt{2}) \sigma $, $\sigma$ the line tension in the corresponding J\"ulicher-Lipowski model \cite{GKRR_MMMAS_2016,EHS_IFB_2022,benes2023degenerate}, $\interfaceparam>0$ is related to the diffuse interface width, and $\dWell(\CHsol) = \frac{1}{4}(\CHsol^2 - 1)^2$ is a double-well potential. 

The second energy we consider accounts for bending properties. As in \cite{fonda2018interface,art:BKV23}, we consider a diffuse interface approximation of the J\"ulicher-Lipowski model \cite{julicher1996shape}:
\begin{equation}
    \label{eq:helfrich}
    \energyHelfr = \int_{\SurfDomain} \frac{1}{2}
    \bendStiff(\CHsol)(\meanCurv-\meanCurv[0](\CHsol))^2\Diff\SurfDomain\,,
\end{equation}
where $\bendStiff(\CHsol)$ is the bending stiffness and $\meanCurv[0](\CHsol)$ the spontaneous curvature. We neglect additional contributions due to the Gaussian curvature of $\SurfDomain$. This is justified as long as the Gaussian bending stiffness is independent on the phase $\CHsol$ and no topological changes occur. To get a continuously differentiable dependency on $\phi$ we consider an interpolation function as in \cite{art:elliott2010surface} and \cite{art:BKV23}: \begin{equation}\label{eq:bendstiff}
  f(\CHsol)=
  \begin{cases}
    f_1 & \mbox{if } \CHsol = 1\\
\displaystyle    \frac{f_1+f_2}{2}
    + \frac{f_1-f_2}{4}\CHsol(3-\CHsol^2) & \mbox{if } -1< \CHsol< 1\\
    f_2 & \mbox{if } \CHsol = -1
  \end{cases}\,,
\end{equation}
for $f \in \{ \bendStiff, \meanCurv[0] \}$ and $f_1,f_2$ the material parameter of $\bendStiff$ or $\meanCurv[0]$ in the separated phases. Together, the energies $\energyGL$ and $\energyHelfr$ define the potential energy $\energy= \energyHelfr+\energyGL$ of the system.

We define the surface material velocity by $\vectvel$ and the kinetic energy as in \cite{reuther2020numerical} by:
\begin{equation}\label{eq:kin}
    \energyKin = \int_{\SurfDomain}
    \frac{\rho}{2}\NORM{\vectvel}^2\Diff\SurfDomain\,,
\end{equation}
with $\rho$ the surface density. For simplicity, we assume $\rho$ to be constant. The Lagrangian  $\Lagrangian=\energyKin-\energy$ is defined as the difference between the kinetic and the potential energy.

We next consider the various sources of dissipation. As in \cite{art:torres-sanchez_millan_arroyo_2019} we define the dissipation potential of the {viscous} stress:
\begin{equation}\label{eq:viscos}
    \DissPotential_V = \int_{\SurfDomain} \viscosity\NORM{\StressTens}^2\,,
\end{equation}
where $\viscosity$ denotes the viscosity and $\StressTens(\vectvel) = \frac{1}{2} (\GradP \vectvel + (\GradP \vectvel)^T)$ is the rate of deformation tensor as considered in \cite{jankuhn2018}. For simplicity, we assume $\viscosity$ to be constant. In addition, we consider the friction with the surrounded material, which is modeled by:
\begin{equation}\label{eq:extResistence}
    \DissPotential_R = \int_{\SurfDomain} \frac{\KinCoef}{2}\NORM{\vectvel}^2\,,
\end{equation}
where $\KinCoef \ge 0$ is a friction coefficient, again assumed to be constant. The third component of the dissipation potential is associated with the dissipation due to phase separation. We assume the immobility potential of the phase field give by:  
\begin{equation}\label{eq:phaseResistence}
    \DissPotential_{\CHsol} = \int_{\SurfDomain} \frac{1}{2\,\mobility} \left\Vert \dot{\CHsol} \right\Vert_{H^{-1}}^2\,,
\end{equation}
with $m>0$ a resistance associated with
$\dot{\phi}=\partial_t\phi+\nabla_{\vectvelV}\phi$ with respect to
$H^{-1}$ norm. Thereby, $\dot{\CHsol}$ denotes the material time
derivative of $\phi$ and $\nabla_{\wb}\phi=(\GradSurf\CHsol,\wb)$, see
Appendix \ref{sec:material_derivative}. The parameter $m$ plays
  the role of a mobility. Here, we follow the approach of
  \cite{MPCMC_JFM_2013} and consider this parameter to be constant. For alternative approaches, we refer to \cite{A_ARMA_2009}. We now define the dissipation potential by $\DissPotential= \DissPotential_V + \DissPotential_R + \DissPotential_{\phi}$. 

Moreover, we add to $\Lagrangian$ and $\DissPotential$ the following constraints. We here only assume local inextensibility of the material as in \cite{art:torres-sanchez_millan_arroyo_2019,reuther2020numerical}. This is incorporated by the Lagrange-function $p$, which serves as the surface pressure and is related to the surface tension. It enters in the constraint by:
\begin{equation}\label{eq:consmass}
  \mathcal{C}_{IE}=-\int_{\SurfDomain}\press\DivP\vectvel\,.
\end{equation}
This constraint induces a conservation of surface area $|\SurfDomain|$. We neglect possible additional constraints, e.g. on the enclosed volume.

With these ingredients, the Lagrange-D'Alembert principle can be
applied. The concept is introduced and used for rigid body models
in~\cite{marsden2001discrete,UDWADIA20021079,IZADI20142570}. Here we
consider this principle for space depended
functions, for which it provides an elegant way to derive
  thermodynamically consistent models. It reads:
\begin{equation}
    0 =  \left( \rho (\partial_t\vectvel+\nabla_{\vectvelV}\vectvel) + \gradFunc{\param}\energy   + \gradFunc{\vectvel}\DissPotential  +  \gradFunc{\vectvel} \mathcal{C}_{IE}, \VecFieldY \right) + \left( \gradFunc{\phi}\energy   + \gradFunc{\dot{\phi}}\DissPotential, \psi \right) + \left( \gradFunc{\press} \mathcal{C}_{IE}, q \right)\,,
\end{equation}
for all test functions $\VecFieldY\in T\REAL^3\vert_{\SurfDomain}$ and $\psi,q\in T^0\SurfDomain$. Here, the symbol $\gradFunc{\cdot}\cdot$ denotes the $L^2$-gradient of the functional. See Appendix~\ref{app:dalembert} for the explicit definitions and calculations. The approach is broadly applicable and provides an alternative to other frameworks, such as the Onsager principle. We obtain the following problem:

\begin{problem}\label{pb:fullmodel-dimensions}
Find $(\vectvel,p,\CHsol,\chempot)$ such that:
\begin{align*}
  \partial_t{\CHsol}+ \nabla_{\vectvelV}\CHsol =& \mobility\,\LapSurf\chempot \,,\\
  \chempot =&\tcurve \left(-\interfaceparam\,\LapSurf\CHsol+\frac{1}{\interfaceparam}\dWell^{\prime}(\CHsol)\right)
  +\frac{1}{2}\bendStiff^{\prime}(\CHsol)\left(\meanCurv-\meanCurv[0](\CHsol)\right)^2 \\
  & -\bendStiff(\CHsol)\meanCurv^{\prime}_0(\CHsol)\left(\meanCurv-\meanCurv[0](\CHsol)\right)\,,\\
  \rho(\partial_t \vectvel + \nabla_{\vectvelV}\vectvel) =& -\GradSurf p -
  p\meanCurv\normalvec + 2\viscosity\DivC\StressTens - \gamma \vectvel + \boldsymbol{b}_{T}+
  \boldsymbol{b}_N \,, {\color{white} \frac{1}{1}}\\
  \DivP\vectvel =& 0 {\color{white} \frac{1}{1}}\,,
\end{align*}
with tangential and normal bending forces defined by:
\begin{align*}
    \boldsymbol{b}_T =& \chempot \GradSurf\CHsol\,, 
    \\
  \boldsymbol{b}_N =& -\left(\LapSurf(\bendStiff(\CHsol)(\meanCurv-\meanCurv[0](\CHsol)))
  +\bendStiff(\CHsol)(\meanCurv-\meanCurv[0](\CHsol))\left(\Vert\shapeOp\Vert^2-\frac{1}{2}\meanCurv(\meanCurv-\meanCurv[0](\CHsol)) \right) \right)
  \normalvec \\ 
  &+\tcurve\left(\frac{\interfaceparam}{2}\Vert\GradSurf\CHsol\Vert^2+\frac{1}{\interfaceparam}
  W(\CHsol)\right)\meanCurv\normalvec   -\tcurve\interfaceparam\GradSurf\CHsol^T\shapeOp\GradSurf\CHsol\normalvec\,, 
\end{align*}
and $\vectvelV$ the relative material
velocity as explained above.
\end{problem}
This problem provides a tight coupling between the phase field $\CHsol$, the surface velocity $\vectvel$, the surface pressure $p$, and the geometry $\SurfDomain$. Exploring the main implications of these couplings is one of the goals of this paper.

\begin{remark}[Thermodynamic consistency]
We refer to Appendix \ref{sec:energy_rate} to show that the relation $$\frac{d}{dt} ( \energyKin + \energy ) = -2 \DissPotential \leq 0 \,,$$ holds for the model in  Problem \ref{pb:fullmodel-dimensions}. This demonstrates thermodynamic consistency. 
\end{remark}

Considering a characteristic length and a characteristic velocity, Problem~\ref{pb:fullmodel-dimensions} can be formulated in non-dimensional form, see Appendix \ref{app:dedim} for details. Keeping the same notation also in non-dimensional form, we obtain:
\begin{problem}\label{pb:fullmodel}
Find $(\vectvel,p,\CHsol,\chempot)$ such that:
\begin{align}
  \partial_t{\CHsol}+ \nabla_{\vectvelV}\CHsol =& \mobility\,\LapSurf\chempot \,,\label{eq:phi}\\
  \chempot=&\tcurve \left(-\interfaceparam\,\LapSurf\CHsol+\frac{1}{\interfaceparam}\dWell^{\prime}(\CHsol) \right)
  +\frac{1}{2}\bendStiff^{\prime}(\CHsol)\left(\meanCurv-\meanCurv[0](\CHsol)\right)^2 \nonumber \\
  &-\bendStiff(\CHsol)\meanCurv^{\prime}_0(\CHsol)\left(\meanCurv-\meanCurv[0](\CHsol)\right)\,,\label{eq:mu}\\
  \partial_t \vectvel + \nabla_{\vectvelV}\vectvel =& -\GradSurf p -
  p\meanCurv\normalvec + \frac{2}{\Reynolds}\DivC\StressTens -\gamma\vectvel+ \boldsymbol{b}_{T}+
  \boldsymbol{b}_N \,,\label{eq:vel}\\
  \DivP\vectvel =& 0 \label{eq:conti}\,,
\end{align}
where $\Reynolds$ denotes the Reynolds number and $\boldsymbol{b}_T$ and $\boldsymbol{b}_N$ are defined as in Problem~\ref{pb:fullmodel-dimensions}.
\end{problem}

In this form, the problem is suitable for numerical approximation. However, before addressing the model in its discrete form, we consider several model simplifications.

\subsection{Model simplifications}
The model in Problem~\ref{pb:fullmodel} is a two component fluid
deformable surface model with phase-dependent elasticity. For
  simplicity, we assume that the material parameters, such as the
  density $\rho$, the viscosity $\viscosity$, and the friction
  coefficient $\gamma$, are constant. It is possible to extend the
  model and consider phase-dependent parameters by following the
  approach described in \cite{AGG_MMMAS_2012}. In the following, we
link our model to simplified models already discussed
in the literature.
\subsubsection*{One component fluid deformable surface} 
In the case of a single phase, the simplified model considers only surface hydrodynamics and its interaction with the geometry due to bending. 
Explicitly, by considering $\CHsol\equiv\mbox{const}$, the system simplifies to:
\begin{align*}
\partial_t \vectvel + \GradSurfConv{\vectvelV}\vectvel &= -\GradSurf p -
  p\meanCurv\normalvec + \frac{2}{\Reynolds}\DivC\StressTens - \gamma\vectvel+
  \boldsymbol{b}_N \,,\\
  \DivP\vectvel &= 0 \,,
\end{align*}
where $\vectvelV=\vectvel-\partial_t\param$ as before, and the normal bending forces reduce to:
\begin{align*}
  \boldsymbol{b}_N =& -\bendStiff\left(\LapSurf\meanCurv
  +(\meanCurv-\meanCurv[0])\left(\Vert\shapeOp\Vert^2-\frac{1}{2}\meanCurv(\meanCurv-\meanCurv[0]) \right) \right)
  \normalvec\,.
\end{align*}
Note that, bending stiffness and spontaneous curvature are now constant parameters. 
In the case of $\meanCurv[0]=0$ and $\gamma=0$, this model has been introduced and simulated in~\cite{reuther2020numerical,art:Krause22} and in the Stokes-limit in \cite{art:torres-sanchez_millan_arroyo_2019}. For further numerical and analytical approaches under additional symmetry assumptions we refer to \cite{AlIzzi2020}, where a linear stability analysis of a tube-geometry is considered, and to \cite{olshanskii2023}, where potential rotational symmetric equilibrium configurations are addressed. For a comparison of different derivations of this model we refer to \cite{Reutheretal_MMS_2015,Reutheretal_MMS_2018} and \cite{brandner2022derivations}.

\subsubsection*{Two component fluid on a stationary surface}

If we assume a stationary surface $\SurfDomain$, i.e. $\vectvel\cdot\normalvec=0$, the model in Problem~\ref{pb:fullmodel} restricts to a pure tangential problem.
We follow the approach of the directional splitting as shown in detail in~\cite{jankuhn2018} and used in~\cite{reuther2020numerical}. We note that $\vectvel_T=\ProjMat\vectvel$ and $u_N=\vectvel\cdot\normalvec$ for the tangential and normal surface velocity, respectively. We project the equations in the surface tangent space and this results in a surface two-phase flow problem:
\begin{align*}
  \partial_t{\CHsol}+ \nabla_{\vectvelV}\CHsol =& \mobility\LapSurf\chempot \,,\\
  \chempot=&\tcurve \left(-\interfaceparam\,\LapSurf\CHsol+\frac{1}{\interfaceparam}\dWell^{\prime}(\CHsol)\right)
  +\frac{1}{2}\bendStiff^{\prime}(\CHsol)\left(\meanCurv-\meanCurv[0](\CHsol)\right)^2 \\
  &-\bendStiff(\CHsol)\meanCurv^{\prime}_0(\CHsol)\left(\meanCurv-\meanCurv[0](\CHsol)\right)\,,\\[0.2em]
  \ProjMat\partial_t \vectvel_T + \nabla_{\vectvelV}\vectvel_T =& -\GradSurf p + \frac{2}{\Reynolds}\divS\StressTens(\vectvel_T) -\gamma\vectvel_T + \mu \GradSurf\CHsol \,,\\
  \divS\vectvel_T =& 0 \,,
\end{align*}
where the relative material velocity simplifies to the Eulerian case $\vectvelV=\vectvel_T$. Because of the tangentiality of the vector and tensor fields, all differential operators fall back to the operators with respect to the covariant derivative. In the case of $\gamma=0$, the model has been introduced and discussed in~\cite{art:BKV23}, where it is used to study the influence of surface hydrodynamics on coarsening in multicomponent scaffolded lipid vesicles (SLV), see \cite{fonda2018interface,art:Fondaetal2019,RFGK_NC_2020}. Without the bending terms, the system has already been addressed in \cite{Nitschke_2012,ABWPK_PRE_2019,Olshanskiietal_VJM_2022}.

\subsubsection*{One component fluid on a stationary surface} 
By considering a single phase on a stationary surface, the system simplifies even further and leads to the inextensible surface Navier-Stokes equations on a stationary surface:
\begin{align*}
\ProjMat\partial_t \vectvel_T + \nabla_{\vectvelV}\vectvel_T &= -\GradSurf p + \frac{2}{\Reynolds}\divS\StressTens(\vectvel_T) - \gamma\vectvel \,,\\
  \divS\vectvel &= 0 \,,
\end{align*}
where again $\vectvelV=\vectvel_T$. These equations have been solved numerically for simply connected surfaces in~\cite{Nitschke_2012,Reutheretal_MMS_2015,Reutheretal_MMS_2018} (and in the Stokes limit in \cite{Olshanskii_Stokes_2018,Brandner2020,Bonito2020}), and for general surfaces in \cite{Nitschke_2017,Reusken2018, Fries2018,reuther2020numerical,Reutheretal_PF_2018,Lederer2020} (and in the Stokes limit in \cite{Olshanskii_Stokes_2019}). 

\subsubsection*{Two component surface without fluid behavior}
A larger literature exists for the evolution of two-component surfaces if surface hydrodynamics is neglected. These models follow as the overdamped limit of Problem \ref{pb:fullmodel}. An explicit derivation is shown in Appendix \ref{app:overdamp}. The resulting model reads:
\begin{align*}
  \begin{aligned} 
  \partial_t{\CHsol}+ \nabla_{\vectvelV}\CHsol =& \tilde\mobility\,\LapSurf\tilde\chempot\,, \\
  \tilde\chempot=&-\tilde\tcurve\interfaceparam\,\LapSurf\CHsol+\frac{\tilde\tcurve}{\interfaceparam}\dWell^{\prime}(\CHsol)\,, \\
  &+\frac{1}{2}\tilde\bendStiff^{\prime}(\CHsol)\left(\meanCurv-\meanCurv[0](\CHsol)\right)^2
  -\tilde\bendStiff(\CHsol)\meanCurv^{\prime}_0(\CHsol)\left(\meanCurv-\meanCurv[0](\CHsol)\right)\,, \\
  u_N =& -\tilde{p}\meanCurv +\tilde{b}_N\,,\\
  \vectvel_T=& -\GradSurf \tilde{p} + \tilde\chempot\GradSurf\CHsol \,,\\
  -\LapSurf\tilde{p} + \tilde{p}\meanCurv^2 & = -\divS \left( \tilde{\chempot}\GradSurf\CHsol\right) + b_N\meanCurv\,,
  \end{aligned}
\end{align*}
where $\tilde{b}_N = \tilde{\boldsymbol{b}}_N\cdot\normalvec$. In this model, $\tilde{p}$ serves as a Lagrange-function to ensure the local inextensibility constraint. The model is similar to a model discussed in \cite{Hauser2013}, see Appendix \ref{app:overdamp} for a detailed comparison. If the constraint on local inextensibility is dropped and only a global area constraint is considered, the model further simplifies. See Appendix \ref{app:overdamp} for the relations with the models considered in  \cite{WD_JMB_2008} and \cite{elliot_stinner2009}.

\section{Numerical discretization}
\label{sec:discrete}

We consider a surface finite element method (SFEM) \cite{art:Dziuk2013,NNV_JCP_2019} to solve the highly nonlinear set of geometric and surface partial differential equations (PDEs) in Problem \ref{pb:fullmodel}. The approach uses higher order surface discretizations, mesh regularization, a Taylor-Hood element for the surface Navier-Stokes equations and established approaches in flat space to split the Navier-Stokes-Cahn-Hilliard like problem. The discretization on the surface addresses numerical analysis results for vector-valued surface PDEs ensuring convergence of surface vector-Laplace and surface Stokes problems on stationary surfaces \cite{hansbo2016,Hardering2022} and reproduces the same expected optimal order of convergences for the full two component fluid deformable surface model.

\subsection{Mesh movement}
By following the work in~\cite{art:Krause22}, we combine the system in Problem~\ref{pb:fullmodel} with a mesh redistribution approach introduced in detail in~\cite{art:Barrett2008} . We consider the initial surface given $\param(0)=\param_0$ and an additional equation for the time evolution of the parametrization:
\begin{align}
        \DerT \param \cdot \normalvec &= \vectvel\cdot\normalvec\,, \label{eq:normalvel}\\
        \meanCurv \normalvec &= \Delta_C \param. \label{eq:meandiff}
\end{align}
This approach generates a tangential mesh movement that maintains the shape regularity and additionally provides an implicit representation of the mean curvature $\meanCurv$. 

\subsection{Surface Approximation}
\label{subsec:surfaceApprox}
We assume that the smooth surface $\SurfDomain$ is approximated by a
discrete $k$-th order approximation $\SurfDomain[\meshparam]$.
Let $\SurfDomain_h^{lin}$ be a piecewise linear reference surface given by shape regular triangulation $\TriangH^{lin}=\{\Cell[i]\}_{i=1}^{\NCell}$. We define $\param$ as bijective map $\param \colon \SurfDomain_h^{lin} \rightarrow \SurfDomain$ such that $\SurfDomain = \cup_{i=1}^{\NCell} \param(\Cell_i)$. The construction of such maps is discussed in \cite{praetorius2020dunecurvedgrid} and \cite{Brandneretal_SIAMJSC_2022}. We get a $k$-th order approximation of $\param$ by the $k$-th order interpolation $\SurfDomain_h^k=\Interpolant^k(\param)$, which defines a higher order triangulation such that $\SurfDomain_h^k = \cup_{i=1}^{\NCell} \param^k_h(\Cell_i)$. We use each geometrical quantity like the normal vector $\normalvec_h$, the shape operator $\shapeOp_h$, and the inner products $(\cdot , \cdot)_h$ with respect the $\SurfDomain_h^k$ where we will drop the index $k$ in the following. We denote the size of the grid by $h$, i.e. the longest edge of the mesh.   

\subsection{Discrete function spaces}

We define the discrete function spaces for scalar function by:
\begin{align}
  V_{k_e}(\SurfDomain_h)=\{ \psi \in C^0(\SurfDomain_h) \vert \psi\vert_{\Cell}\in\mathcal{P}_{k_e}(\Cell)\}\,,
\end{align}
with $\mathcal{P}_{k_e}$ the space polynomials, where we set the element order as $k_e=k$. We define $\boldsymbol{V}_{k}(\SurfDomain[\meshparam])=[V_{k}(\SurfDomain[\meshparam])]^3$ as space of discrete vector fields. We discretize $\vectvel_h,\param_h\in\boldsymbol{V}_2(\SurfDomain_h)$,  $\meanCurv_h,\CHsol_h,\chempot_h\in V_2(\SurfDomain_h)$, and $p_h\in V_1(\SurfDomain_h)$. This corresponds to a Taylor-Hood element for the pair of velocity and pressure. 

\subsection{Discrete model}
Following the strategy in \cite{art:BKV23}, we separate the surface phase field and surface Navier–Stokes equations by an operator splitting approach. In contrast with \cite{art:BKV23}, where the surface was stationary, here the surface evolves and we adapt the numerical schemes used in \cite{art:Krause22} for one-component systems.

Let $\{ t^n \}_i^N$ be time interval discretization, with $t^n=n\tau$ and $\tau>0$ be the time steps. We first solve the surface phase field problem, eqs.~\eqref{eq:phi}-\eqref{eq:mu}, and then solve the surface Navier-Stokes, eqs. \eqref{eq:vel} and \eqref{eq:conti}, together with eqs. \eqref{eq:normalvel}-\eqref{eq:meandiff} for mesh regularization.
The inner product $\InnerApprox{\cdot}{\cdot}$ is approximated by a quadrature rule with an order that is chosen high enough such that test and trial
functions and area elements are well integrated. We denote the time discrete surface by $\SurfDomain_h^{n-1}=\SurfDomain_h(t^{n-1})$. 

The equations are discretized in time by a semi-implicit Euler scheme where the nonlinear terms are chosen explicitly apart from the double well potential. As is common for Cahn-Hilliard equations, we linearize the derivative of the double-well potential by a Taylor expansion of order one.

\begin{problem}[Discrete surface Cahn-Hilliard problem] \label{pb:CHdis}
Find $(\CHsolApprox,\chempotApprox)\in[\TestSpaceApprox\times\TestSpaceApprox](\SurfDomain_h^{n-1})$ such that:
\begin{align*}
\frac{1}{\tau}\InnerApprox{\CHsolApprox^n}{\TestCHApprox}
+ \InnerApprox{\GradSurfConv{\vectvelV^{n-1}}\CHsolApprox^n}{\TestCHApprox}
=& \frac{1}{\tau}\InnerApprox{\CHsolApprox^{n-1}}{\TestCHApprox}
- \mobility \InnerApprox{\GradSurf\chempotApprox^{n}}{\GradSurf\TestCHApprox} \,,
\\[0.7 em]
\InnerApprox{\chempotApprox^n}{\TestCHmuApprox}
=& \, \interfaceparam \InnerApprox{\GradSurf\CHsolApprox^n}{\GradSurf\TestCHmuApprox} 
+ \frac{1}{\interfaceparam} \InnerApprox{-2(\CHsolApprox^{n-1})^3 +
\left(3(\CHsolApprox^{n-1})^2-1\right)\CHsolApprox^n }{\TestCHmuApprox}\\
&
+ \frac{1}{2}\InnerApprox{
   \bendStiff^{\prime}(\CHsolApprox^{n-1})\left(\meanCurvApprox^{n-1}-\meanCurv[0](\CHsol)(\CHsolApprox^{n-1})\right)^2}{\TestCHmuApprox} \\
& 
- \InnerApprox{\bendStiff(\CHsolApprox^{n-1})\meanCurv[0](\CHsol)^{\prime}(\CHsolApprox^{n-1})\left(\meanCurvApprox^{n-1}-\meanCurv[0](\CHsol)(\CHsolApprox^{n-1})\right)}{\TestCHmuApprox} \,,  
\end{align*}
for all $(\TestCHApprox,\TestCHmuApprox)\in [\TestSpaceApprox\times\TestSpaceApprox](\SurfDomain_h^{n-1})$.
\end{problem}

The surface Navier-Stokes equations and the mesh redistribution are considered together. Moreover, we define a discrete surface update variable $\update_h^{n}=\param_h^{n}-\param_h^{n-1}$, 
which is considered as unknown instead of the surface parametrization $\param^{n}$.
For time discretization, we again consider a semi-implicit Euler scheme
following the strategy presented
in~\cite{art:Barrett2008,art:Krause22}.
We obtain the following discrete problem.
\begin{problem}[Discrete surface Navier-Stokes and surface update problem] \label{pb:NSdis}
Find $(\vectvelApprox,\pressApprox,\meanCurv_h,\updateApprox)\in[\TestSpaceVec\times\TestSpaceApprox\times\TestSpaceApprox\times\TestSpaceVec](\SurfDomain_h^{n-1})$ such that: 
\begin{align*}
  \frac{1}{\tau} \InnerApprox{\vectvelApprox^{n}}{\TestVecApprox} + \InnerApprox{\GradSurfConv{\vectvelV_{\meshparam}^{n-1}}\vectvelApprox^n}{\TestVecApprox} = &
  \InnerApprox{\pressApprox^n}{\DivP\TestVecApprox} 
   - \frac{2}{\Reynolds} \InnerApprox{
   (\StressTens(\vectvelApprox^{n})}{\GradP\TestVecApprox} \\
   & + \InnerApprox{\chempotApprox^n\GradSurf\CHsolApprox^n}{\TestVecApprox} \\ 
   & + \InnerApprox{ \GradSurf\left(\kappa(\CHsol^{n})(\meanCurv_h-\meanCurv_0(\CHsol^{n}))\right) }{\GradSurf(\TestVecApprox\cdot\normalvec_h)} \\
   &-\InnerApprox{\kappa(\CHsol^{n})(\meanCurv_h-\meanCurv_0(\CHsol^{n}))B^{n-1}}{\TestVecApprox\cdot\normalvec_h} \\
  &+\frac{1}{\tau} \InnerApprox{\vectvelApprox^{n-1}}{\TestVecApprox}
  -\InnerApprox{\gamma\vectvel_h^n}{\TestVecApprox}
 \,,
 \\
\InnerApprox{\DivP\vectvelApprox^{n}}{\TestPressApprox} = & \;0 \,,\\[1em]
    \frac{1}{\tau}\InnerApprox{\update_h^n \cdot \normalvecApprox}{\TestHApprox} =& \,\InnerApprox{\vectvel_h^n\cdot\normalvecApprox}{\TestHApprox}  \,,\\
    \InnerApprox{\meanCurv_h^n\normalvec_h}{\TestZApprox} 
    + \InnerApprox{\bendStiff(\CHsol_h^n)\GradC \update_h^n}{\GradC \TestZApprox} =& 
    -\InnerApprox{\bendStiff(\CHsol_h^n)\GradC \MapU_h^{n-1}}{\GradC \TestZApprox}\,, \\
\end{align*}
for all $(\TestVecApprox,\TestPressApprox,\TestHApprox,Z_{\meshparam})\in[\TestSpaceVec\times\TestSpaceApprox\times\TestSpaceApprox\times\TestSpaceVec](\SurfDomain_h^{n-1})$, where \\
$B^{n-1}=\left(\Vert\shapeOp_h\Vert^2-\frac{1}{2}\tr \shapeOp_h(\tr\shapeOp_h-\meanCurv[0](\phi^n)) \right)$.\footnote{We note a typing error in \cite{art:Krause22}, which has been confirmed by the authors.}
\end{problem}

The solutions of the discrete problems above are intermediate solutions with respect to the old surface $\SurfDomain^{n-1}_h$. For each variable $\hat\psi\in V_h(\SurfDomain_h^{n-1})$, we define the lift $\psi\in V_h(\SurfDomain_h^{n})$ by the evaluation with respect to the linear reference geometry where $\psi(\param_h^n) = \hat\psi(\param_h^{n-1})$. This corresponds to a nodal interpolation with respect to the degrees of freedom. We consider the following steps in each time step:
\begin{enumerate}
    \item Compute $(\hat\CHsol_h,\hat\chempot_h)\in[\TestSpaceApprox\times\TestSpaceApprox](\SurfDomain_h^{n-1})$ as intermediate solution of Problem \ref{pb:CHdis}.
    \item Compute $(\hat\vectvel_h,\hat\press,\hat\meanCurv_h,\hat\update_h)\in[\TestSpaceVec\times\TestSpaceApprox\times\TestSpaceApprox\times\TestSpaceVec](\SurfDomain_h^{n-1})$ as intermediate solution of Problem \ref{pb:NSdis}. 
    \item Compute the new surface $\SurfDomain^{n}_h$ by updating its parametrization $\param_h^n=\param_h^{n-1}+\update_h^{n}$.
    \item Lift the intermediate solutions on the new surface to get $(\CHsol_h,\chempot_h)\in[\TestSpaceApprox\times\TestSpaceApprox](\SurfDomain_h^{n})$ and $(\vectvel_h,\press,\meanCurv_h,\update_h)\in[\TestSpaceVec\times\TestSpaceApprox\times\TestSpaceApprox\times\TestSpaceVec](\SurfDomain_h^{n})$.
\end{enumerate}

\subsection{Implementational aspects}
The discrete systems are implemented within the finite element toolbox AMDiS \cite{vey2007amdis,witkowski2015software} based on DUNE \cite{SanderDune2020,alkamper2014dune}. The surface approximation is done by the Dune-CurvedGrid library developed in \cite{praetorius2020dunecurvedgrid}. For a straightforward mesh parallelization and  multiple processor computation we used the PETSc library and solved the linear system by using a direct solver. Due to the nature of the equations, which lack a maximum principle, we allow $\phi_h \in \REAL$ and extend all material parameters constantly for $\phi_h < -1$ and $\phi_h>1$.

The SFEM approach discussed in Section \ref{subsec:surfaceApprox} does not support topological changes of the domain. Such changes can occur due to pinch offs associated with budding events or the formation of furrows. The pinch offs are characterized by a negative Gaussian curvature $\GaussCurv$. We use this property as indicator to define our simulated time interval $[0,T]$, where the final simulation time $T$ is defined by $T=\inf\{ t>0\vert \GaussCurv(t) < \GaussCurv_0\}$ where $\GaussCurv_0<0$ is a chosen lower bound. If this criteria is not met, we consider $T=\infty$ and solve until an equilibrium configuration is reached. Numerical methods which are able to handle topological changes, or at least able to recover after such an event in a meaningful physical state, require an implicit description of the surface. Candidates are trace finite element methods (TraceFEM) \cite{10.1093/imanum/drz062}, cut finite element methods (CutFEM), level-set methods or diffuse interface methods \cite{nestler2018orientational}. However, these methods are computationally more expensive and currently not applicable for the considered problem. A comparison between SFEM and TraceFEM for Stokes flow on stationary surfaces shows a factor of $10^2 - 10^3$ lower errors for SFEM for comparable mesh sizes, depending on the considered error measure \cite{Brandneretal_SIAMJSC_2022}. Applications of CutFEM and level set methods for this problem class are not known and the first detailed numerical investigations for diffuse interface methods for vector-valued surface PDEs show strong limitations for the required accuracy for the reconstructed geometric terms \cite{nestler2023diffuse}. These results are supported in a benchmark problem for vector-valued diffusion on surfaces with large curvature \cite{bachini2022diffusion}. This motivates the use of SFEM and the needs to consider simulations only before a potential topological change.

\section{Numerical results}
\label{sec:results}

\subsection{Considered parameters and initial conditions}

Let $\SurfDomain$ be the unit sphere with radius $R = 1$. The velocity field is initialized by $\vectvel_0=0$ and we define two different initial configurations for the phase field. The first one is $\CHsol_0(\coord)=\tanh(\alpha x_0)$ that defines a symmetric and equal distributed phase field where the phases are separated. The second one defines a random initial condition: $\phi_1 = \min\{\max \{\hat{\phi}_1,-1\},1\}$, where $\hat{\phi}$ is given by Gaussian random field:
\begin{align*}
    \hat{\phi}_1(\coord)=-1.0 + \sum\limits_{i=0}^N e^{-\frac{\beta}{2}\Vert \coord - \coord_i\Vert^2} \quad \forall \coord\in\SurfDomain\,,
\end{align*}
where $\{\coord_i\}_0^N \subset \SurfDomain$ and $\alpha=\beta=100$. This allows to create a reproducible random initial condition. Both configurations consider an equal distributed phase field. We vary the Reynolds number $\Reynolds$ and the bending stiffness $\kappa$ but set $\meanCurv_0=\gamma=0$. Other parameters are set to $\epsilon = 0.02$, $m = 0.001$ and $\sigma=1.0$ ($\tilde{\sigma} = \frac{3}{2} \sqrt{2}$), which provides a good compromise between computational effort and physical accuracy. Numerical parameters are such that $h^3 \sim \tau$.

\subsection{Convergence study}

Due to the tight coupling between the phase field $\CHsol$, the surface velocity $\vectvel$, the surface pressure $p$, and the geometry $\SurfDomain$ we expect the numerical solution to sensitively depend on the approximations. Therefore, we start by considering a convergence study. Due to a lack of analytical solutions for the full problem, we compute a numerical solution with a fine discretization and study convergence with respect to this solution. In addition, we address general properties of the solution.

\begin{figure}
    \centering
    \includegraphics[width=\textwidth]{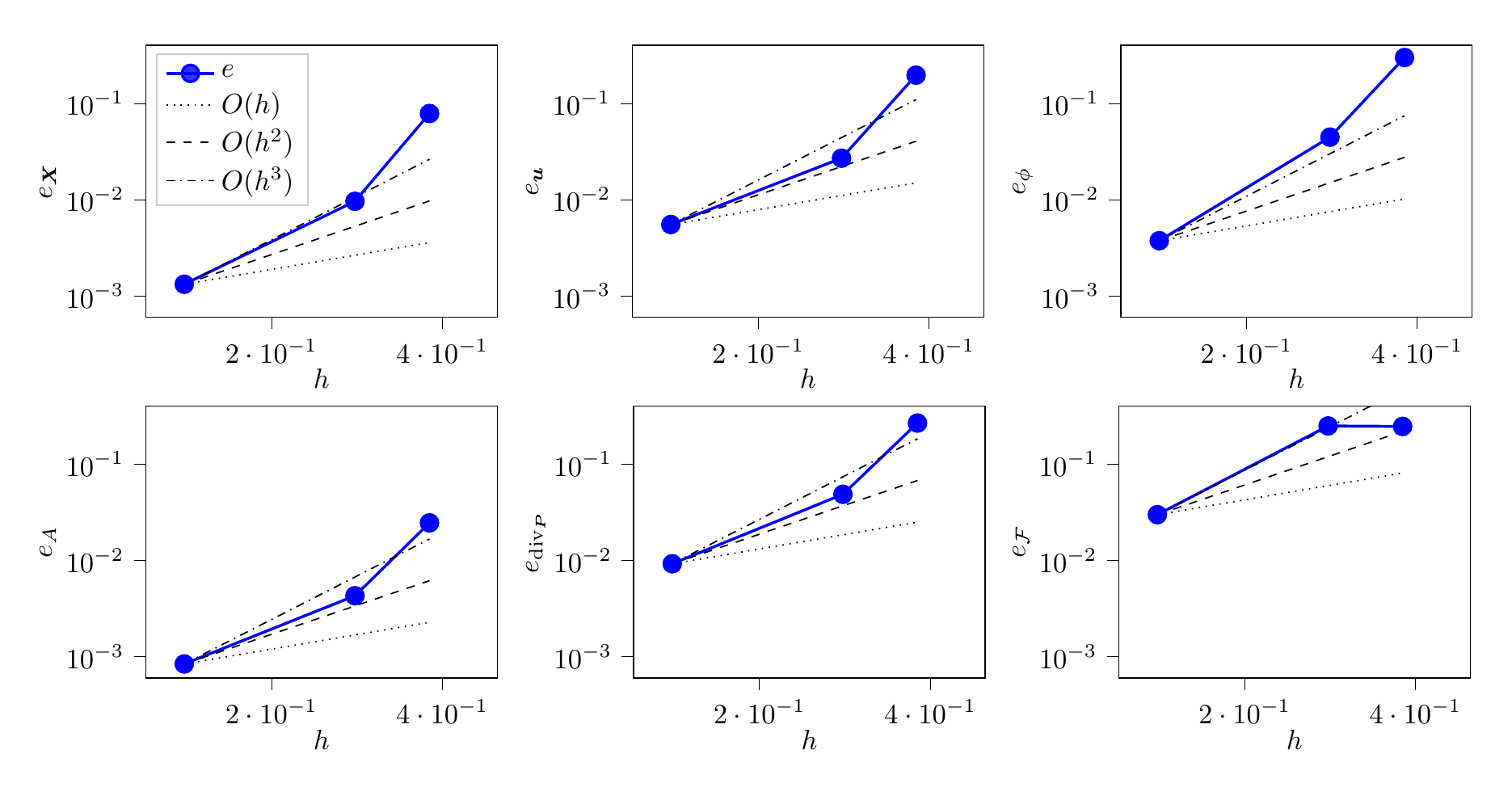}
    \caption{Convergence study for two-phase fluid deformable surfaces with respect to the mesh size $h^3 \sim \tau$ for $\param_h,\vectvel_h,\phi_h$ (top row) and error for area $A$, divergence $div_\mathbf{P} \mathbf{u}_h$ and total energy ${\cal{F}}_K + {\cal{F}}$ (bottom row).}
    \label{fig:convergenceStudy}
\end{figure}

We compute the errors with respect to the $L^2$ norm in space and the $L^{\infty}$ norm in time. The norms are calculated with respect to the linear reference geometry $\SurfDomain_h^{lin}$, where the evaluation of the norms on different surfaces $\SurfDomain,\SurfDomain_h,\SurfDomain_h^{lin}$ leads to equivalent norms with the same order of convergence \cite{art:Demlow2009}. We define the following errors:
\begin{align*}
  e_{\param}&=\Vert \param_h - \param \Vert_{L^{\infty}(L^2(\SurfDomain_h^{lin}))} , \\
  e_{\vectvel}&=\Vert \vectvel_h(\param_h) - \vectvel(\param) \Vert_{L^{\infty}(L^2(\SurfDomain_h^{lin}))} , \\
  e_{\phi}&=\Vert \CHsol_h(\param_h) - \CHsol(\param) \Vert_{L^{\infty}(L^2(\SurfDomain_h^{lin}))}, \\
  e_A &= \vert \Vert\SurfDomain_h(t)\vert - \vert \SurfDomain_h(0) \vert \vert_{L^{\infty}}, \\
  e_{\DivP}&=\Vert \DivP\vectvel_h \Vert_{L^{\infty}(L^2(\SurfDomain_h))} , \\  
  e_{\energy}&=\vert (\energyKin_h+\energy_h) - (\energyKin+\energy) \vert_{L^{\infty}} , 
\end{align*}
with discrete solutions $\param_h,\vectvel_h,\phi_h$ with varying $h$ and $\param,\vectvel,\phi$ the discrete reference solutions with the finest mesh size. 
Additionally, we measure the area conservation error $e_A$, the local inextensibility error $e_{\DivP}$, and the error of the total energy $e_{\energy}$, where again $\energyKin$ and $\energy$ denote the discrete reference solution on the finest mesh size. 
The simulations are done for the random initial condition $\CHsol_1$ with $\Reynolds = 1.0$ and $\kappa_1 = \kappa_2 = \kappa = 0.02$. 

The results are shown in Figure \ref{fig:convergenceStudy}. They indicate third order convergence for the surface error $e_{\param}$ and the error of the phase field $e_{\CHsol}$, which corresponds to the results in \cite{praetorius2020dunecurvedgrid,art:Demlow2009} for simple mean curvature flow and scalar-valued surface PDEs. The order of convergence for the phase field $e_{\CHsol}$ also agrees with the one obtained for the corresponding problem on stationary surfaces considered in \cite{art:BKV23}. For the error of the velocity field $e_{\vectvel}$ the results indicate second order convergence in accordance with the results in \cite{art:Krause22} for one-component fluid deformable surfaces. In \cite{Brandneretal_SIAMJSC_2022}, third order convergence of $e_{\vectvel}$ is shown for the tangential flow of the surface Stokes equations on stationary surfaces. However, this requires a higher order approximation of the normal vector \cite{hansbo2016,Hardering2022}, which is not fulfilled in our approach. Furthermore, it remains open if this increased order also emerges for the surface Navier-Stokes equations and on evolving surfaces. For the area conservation error $e_A$ and the inextensibility error $e_{\DivP}$ we get second order convergence with respect to $h$. Both are related to each other. While $e_A = 0$, the same order for $e_{\DivP}$ has been obtained in \cite{art:BKV23} for a stationary surface. The convergence error of the total energy indicates third order. This might be due to the dominating effect of Ginzburg-Landau energy as the expected order for the underlying Helfrich model would be lower \cite{Dziuk_NM_2008}. However, overall we see  experimentally the expected optimal order of convergence for the full problem.

\subsection{Phase separation, bulging and induced flow}
\begin{figure}
    \centering
    \includegraphics[width=\textwidth]{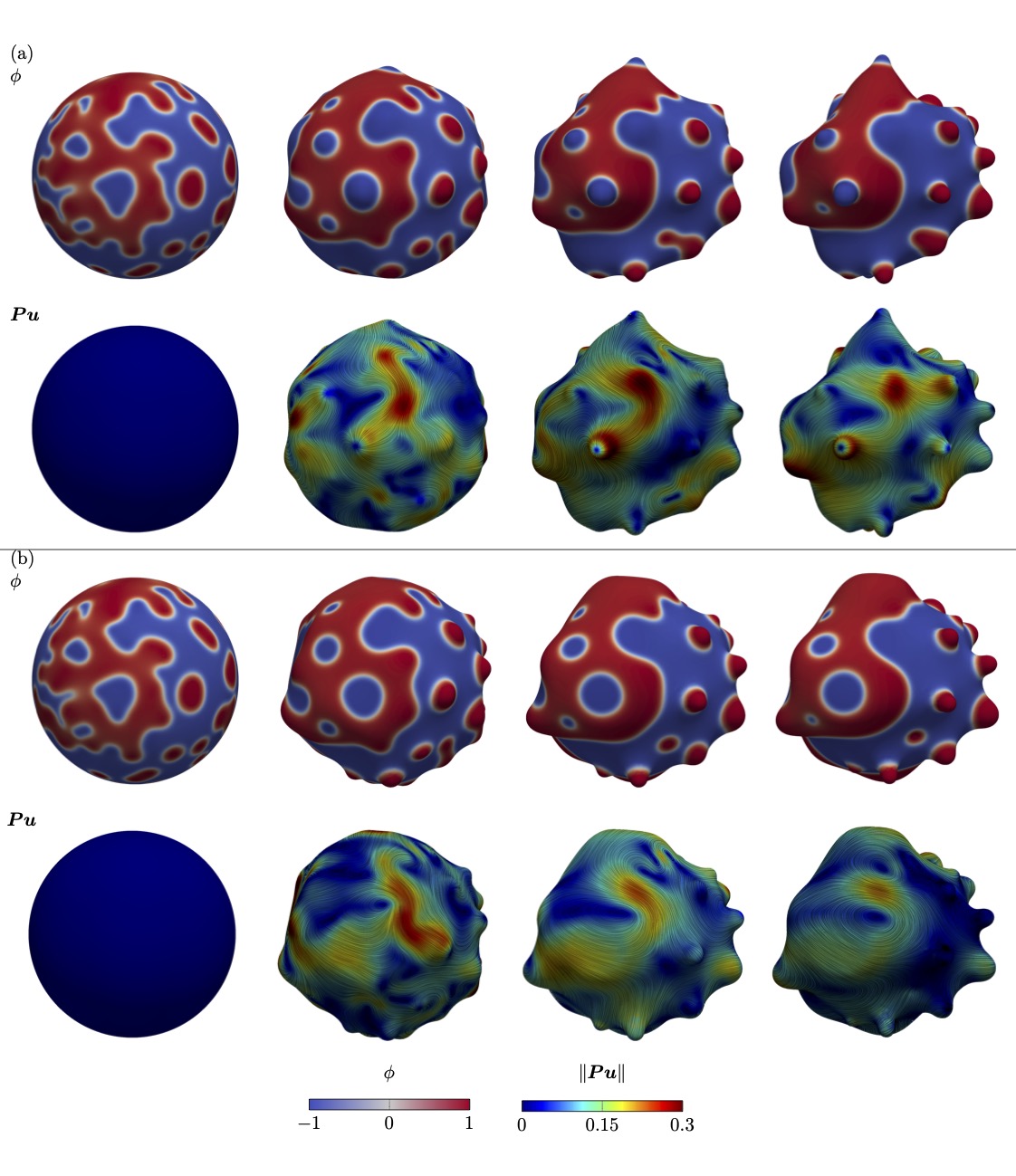}    
    \caption{(a),(b) Snapshots of the relaxation of the two-component fluid deformable surface with random initial condition $\CHsol_1$ and $\Reynolds=1.0$ for $t=0,0.3,0.8,1.1$ (from left to right), with constant bending stiffness $\bendStiff_1 = \bendStiff_2 = \bendStiff=0.02$ in (a) and phase-depended bending stiffness $\bendStiff_1 = 0.02$ and $\bendStiff_2 = 0.5$ in (b). In (b) the red colored phase is less stiff than the blue colored phase. Top row: phase field $\phi$. Bottom row: tangential velocity $\ProjMat\vectvel$. The flow is visualized by a LIC filter and color coding represents the magnitude. Corresponding movies are provided in the supplementary data. 
    }  
    \label{fig:evolution}
\end{figure}
\begin{figure}
    \centering
    \includegraphics[width=\textwidth]{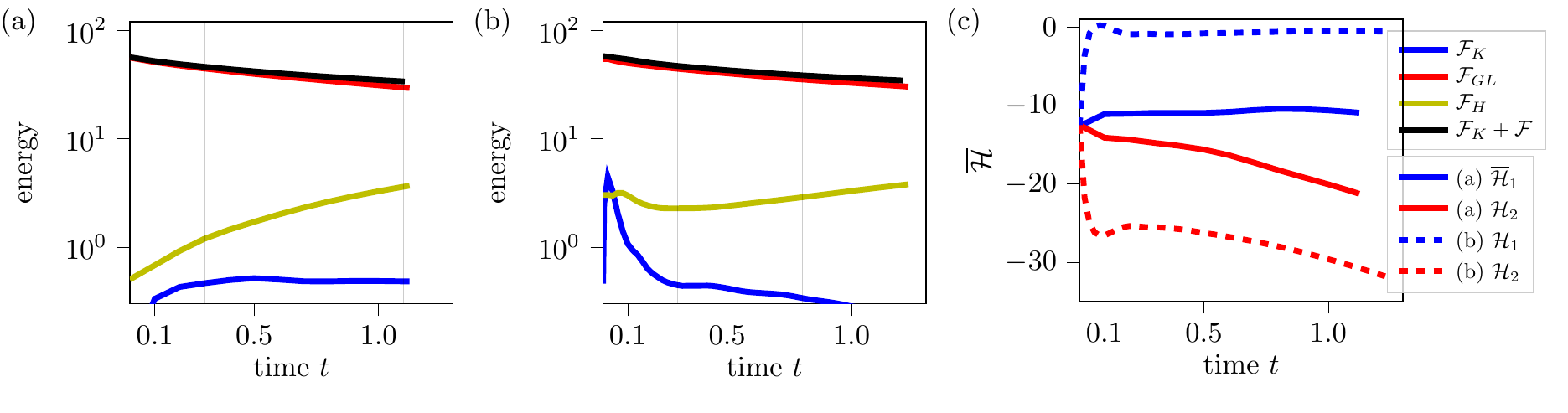}    
    \caption{(a),(b) The energies $\energyKin$, $\energyGL$, $\energyHelfr$ and $\energy_K+\energy$ over time where (a) corresponds to fig \ref{fig:evolution} (a) and (b) corresponds to fig \ref{fig:evolution}  (b). The time instances are highlighted in the plots. (c) Averaged mean curvature $\overline{\meanCurv}_{1,2}$ for the different phases with respect to the simulations done in fig \ref{fig:evolution} (a) and (b), the color corresponds to the colored phases.
    }  
    \label{fig:energies}
\end{figure}
To demonstrate the strong interplay between composition, curvature, and hydrodynamics we consider $\phi_1$ as initial condition. The other parameters are $\Reynolds = 1.0$ and either $\bendStiff_1 = \bendStiff_2 = \bendStiff = 0.02$ or $\bendStiff_1 = 0.02$ and $\bendStiff_2 = 0.5$. In the last case, the red colored phase has a lower bending stiffness and is therefore expected to be guided to or initiate regions with higher mean curvature, while the blue color phase has a larger bending stiffness and can be expected to prefer regions of lower mean curvature. The evolution is shown in Figure \ref{fig:evolution} a) and b), respectively. Both cases show the composition and the tangential velocity for selected time instances. In both cases red and blue islands form and a wavy interface between larger red and blue regions is established. Some coarsening events can be spotted and the wavy interface flattens over time. However, the main evolution is in the shape, which strongly deforms and forms bulges. Especially for circular islands the interface length is reduced by bulging leading to strong deformation from the sphere. These shape changes, but also the coarsening events induce flow. The simulations are only shown for a relatively short period of time and are terminated before a potential topological change would happen. Differences between a) and b) are also visible. While in a) the red and blue phases evolve similarly, in b) the blue phase, the one associated with a larger bending stiffness, forms less curved regions. Especially islands of this phase do not bulge out. Instead, they become relatively flat. The behavior of the simulation shown in Figure \ref{fig:evolution} are quantified in Figure \ref{fig:energies} by different energies. In the current setting, the evolution is dominated by $\energyGL$, which is reduced by coarsening but also by bulging, which is associated with an increase in $\energyHelfr$. Hydrodynamics seems to play a minor role in the evolution. $\energyKin$ is 1-2 orders of magnitude lower than $\energyGL$. The strong increase of $\energyKin$ and $\energyHelfr$ at the beginning is associated with the increased driving force resulting from the phase-dependent bending stiffness. The spatial average over the mean curvature for the red and the blue phase, i.e.,
\begin{align*}
    \overline{\meanCurv}_1 = \int_{\{\CHsol<0\}} \meanCurv 
    \quad \text{and} \quad
    \overline{\meanCurv}_2 = \int_{\{\CHsol>0\}} \meanCurv\,,
\end{align*}
is shown in Figure 3 c). The deviation between $\overline{\meanCurv}_1$ and $\overline{\meanCurv}_2$ for the setting of Figure \ref{fig:evolution} a) results from asymmetries in the initial condition. For the setting of Figure \ref{fig:evolution} b), these quantities show strong differences between the different phases. The blue phase on average forms flat regions, while the red phase forms strongly curved bulges. 

\subsection{Variation of parameters}
\begin{figure}
    \centering
    \includegraphics[width=\textwidth]{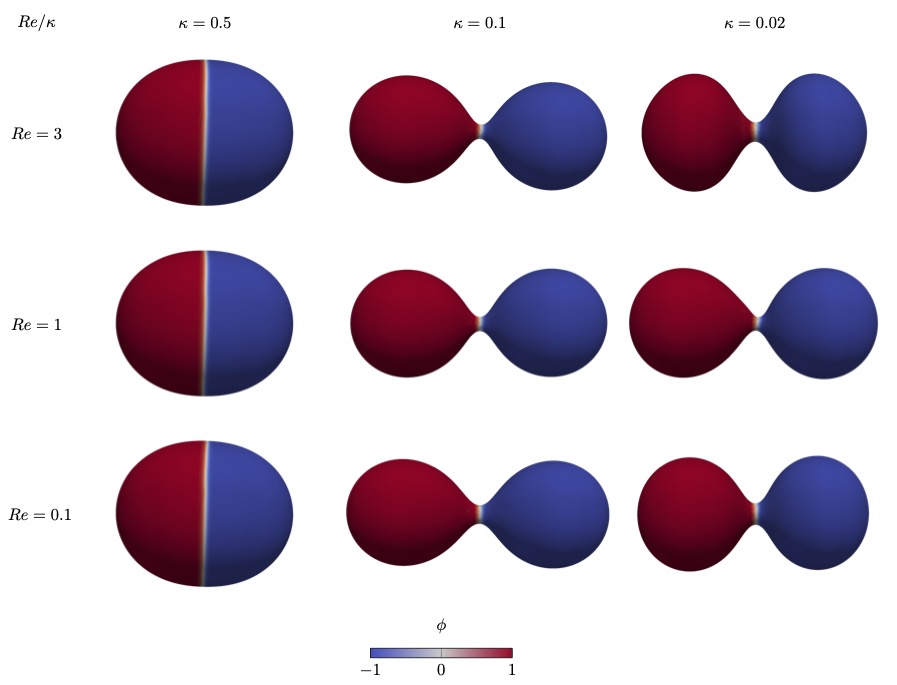}
    \caption{Final configuration obtained for different parameters and initial condition $\CHsol_0$. Either the simulation reaches the equilibrium configuration or the criteria for a potential pinch-off is reached. The corresponding critical times $T$ are shown in Table \ref{tab:criticalTimes}.}
    \label{fig:relaxBean}
\end{figure}

\begin{table}
  \centering
  \begin{tabular}{c c| c|c|c}
    (a)&$\Reynolds$/$\kappa$  & \textbf{0.5} & \textbf{0.1} & \textbf{0.02} \\
    &\textbf{3} & $\infty$ & 5.6   & 1.9      \\
    &\textbf{1} & $\infty$ & 10.7  & 5.3      \\
    &\textbf{0.1} & $\infty$ & 87.0  & 42.0
  \end{tabular}  
  \begin{tabular}{c c| c|c|c}
    (b)&$\Reynolds$/$\kappa$  & \textbf{0.5} & \textbf{0.1} & \textbf{0.02} \\
    &\textbf{3} & $\infty$ & 9.1 & 0.5      \\
    &\textbf{1} & $\infty$ & 20.0& 1.1      \\
    &\textbf{0.1} & $\infty$ & 170.5 & 9.8
  \end{tabular}  
  \hspace{1cm}
  \caption{Critical times $T$ for the different Reynolds numbers $\Reynolds$ and bending stiffness' $\bendStiff$ for the symmetric initial value $\CHsol_0$ in (a) and the random initial value $\phi_1$ in (b).}
  \label{tab:criticalTimes}
\end{table}

\begin{figure}
    \centering
    \includegraphics[width=\textwidth]{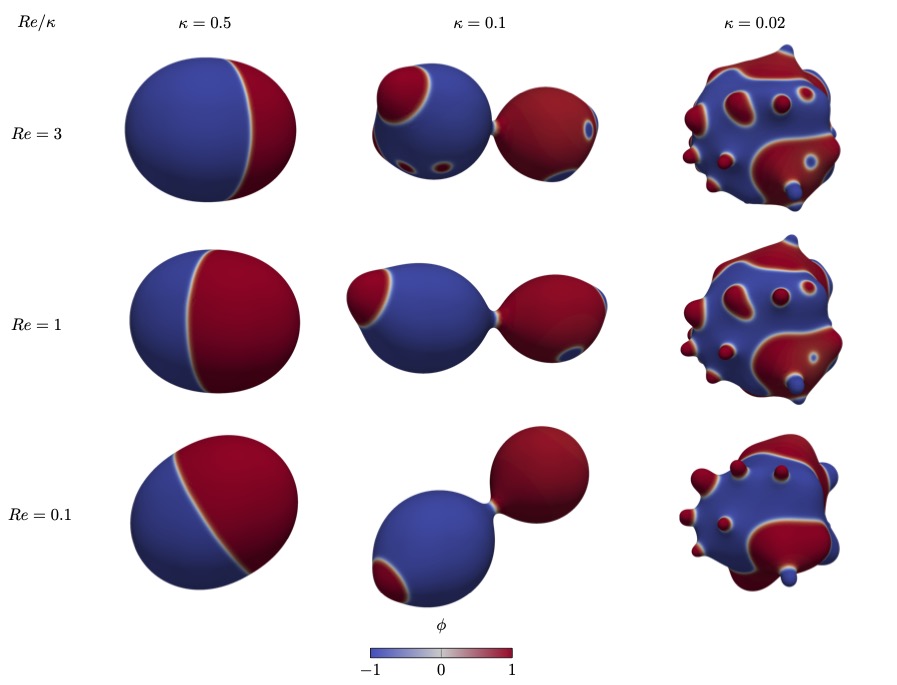}
    \caption{Final configuration obtained for different parameters and initial condition $\CHsol_1$. Either the simulation reaches the equilibrium configuration or the criteria for a potential pinch-off is reached. The corresponding critical times $T$ are shown in Table \ref{tab:criticalTimes}.}
    \label{fig:relaxRandom}
\end{figure}

We now explore the influence of the parameters in more details. To reduce the number of parameters we only consider the situation of a constant bending stiffness $\bendStiff_1 = \bendStiff_2 = \bendStiff$, set $\bendStiff = \{0.02, 0.1, 0.5\}$, and vary the Reynolds number with $\Reynolds = \{0.1, 1, 3\}$. We consider both initial conditions $\CHsol_0$ and $\CHsol_1$. The behavior is demonstrated for the symmetric initial phase field $\CHsol_0$ in Figure \ref{fig:relaxBean} and for the random initial phase field $\CHsol_1$ in Figure \ref{fig:relaxRandom}. We only show the final configuration, this is either the equilibrium configuration or the configuration at the critical time $T$ before a potential pinch-off. A table with critical final times is presented in Table \ref{tab:criticalTimes}. In the first case ($\CHsol_0$) the final configuration is mainly determined by the interplay of the Helfrich energy $\energyHelfr$ and the Ginsburg-Landau energy $\energyGL$. For $\bendStiff = 0.5$, $\energyHelfr$ dominates. Any deformation of the initial sphere requires to increase the Helfrich energy. Only small deformations are possible. The interface length, and thus $\energyGL$, is reduced by deforming the sphere to an ellipsoid-like shape with the interface placed along the short axis. Further reduction of the interface length either requires to increase the mean curvature at the tips of the long axis of this shape or to form a furrow alone the interface. Both effects strongly increase $\energyHelfr$. This is only realized for settings with lower $\bendStiff$ and eventually will lead to a pinch-off. Hydrodynamics seems to have no significant effect in these evolutions.
However, a closer look at the consider critical times $T$ in Table \ref{tab:criticalTimes} shows the opposite. While the final configurations look similar, the time to reach these configurations is strongly influenced by hydrodynamics. For low bending stiffness, $\bendStiff = \{0.1, 0.02\}$, the critical time is drastically reduced if the Reynolds number $\Reynolds$ is increased. 

For the random initial condition $\phi_1$ the results are shown in Figure \ref{fig:relaxRandom}. Here, we have an interplay between coarsening, shape deformation, bulging, and furrow formation. For large bending stiffness $\bendStiff = 0.5$, where strong shape deformations are suppressed, the same final configuration as in Figure \ref{fig:relaxBean} is reached, meaning an ellipsoidal like shape with the interface placed along the short axis. This changes for lower $\bendStiff$. For $\bendStiff=0.1$, we observe partial coarsening, bulging of islands, and the  formation of a furrow, which eventually leads to a pinch-off. The configurations also change with hydrodynamics. For $\Reynolds = 3$ the potential pinch-off happens earlier. Several islands are still present in the red as well as the blue phase. Decreasing $\Reynolds$ also leads to possible pinch-offs but at a later coarsening state. There are less islands present. While it is known in flat space that hydrodynamics can enhance coarsening and this is also demonstrated on stationary surfaces \cite{art:BKV23}, these results indicate that this also holds on evolving surfaces and that hydrodynamics enhances furrow formation and potential pinch-off. This is confirmed in Table \ref{tab:criticalTimes}, which again shows a drastic reduction of the critical time $T$. For $\bendStiff = 0.02$ the evolution is dominated by bulging. As discussed in the previous section bulging is associated with large local absolute mean curvature values but allows to reduce the interface length. The critical time results from potential pinch-offs of the bulges. Coarsening is only a minor effect. Again hydrodynamic drastically enhances the evolution.  
Increasing $\Reynolds$ again drastically reduces the time for potential pinch-off, see Table \ref{tab:criticalTimes}. 

\section{Conclusions}
\label{sec:conclusion}

While the literature on coexisting fluid domains in model systems for biomembranes and their dynamics is rich, the influence of surface viscosity has not been discussed in this context. We have filled this gap and derived a thermodynamically consistent model that accounts for surface hydrodynamics. The derivation of the model by a Lagrange-D'Alembert principle is applicable in general. Simplifications of the derived two-phase fluid deformable surface model lead to known models in the literature, and this further confirms the validity of the approach.

By combining numerical approaches for fluid deformable surfaces \cite{art:Krause22} and two-phase flow models on stationary surfaces \cite{art:BKV23}, we obtain a numerical approach for the full model, which is demonstrated to converge with expected optimal order. This provides the basis for a detailed investigation of the strong interplay of surface phase composition, surface curvature, and surface hydrodynamics. Depending on the material parameters, the line tension $\sigma$, the Reynolds number $\Reynolds$, and the bending stiffness $\bendStiff$ the evolution is dominated by coarsening or shape evolution. Both phenomena are shown to be strongly enhanced by hydrodynamics. In situations where the line tension and the bending stiffness are compatible, the interface length is reduced by the formation of a furrow or the formation of bulges. Both potentially lead to pinch-offs. The enhanced evolution towards such topological changes driven by hydrodynamics has various biological implications. In the context of a furrow formation and its subsequent shrinkage, this can be associated with cell division \cite{Mietkeetal_PRL_2019}; in the context of bulges with pinch-offs, this can be associated with endocytosis and exocytosis \cite{Kasonenetal_NRMCB_2018,AlIzzi2020}. While quantitative comparisons with experimental data in these context require further studies, the mathematical model and the numerical algorithm to solve it with the required accuracy are provided.

\paragraph{Acknowledgements} We acknowledge computing resources provided by ZIH at TU Dresden and by JSC at FZ J\"ulich, within projects WIR and PFAMDIS, respectively.

\paragraph{Funding} This work was supported by the German Research Foundation (DFG) within the Research Unit ``Vector- and Tensor-Valued Surface PDEs'' (FOR3013).

\paragraph{Declaration of interests} The authors report no conflict of interest.

\paragraph{Data availability statement} The data that support the findings of this study are available from the corresponding author upon reasonable request.

\paragraph{Author ORCIDs}E. Bachini, https://orcid.org/0000-0001-8610-9426;\\
V. Krause, https://orcid.org/0009-0006-3437-5902 \\
A. Voigt, https://orcid.org/0000-0003-2564-3697


\appendix

\section{Symbols, definitions, and geometric notation}
\label{app:symbols}

We define the surface differential operators as used in Section \ref{sec:model}, pointing out the differences if they are applied to scalar, vector, or tensor quantities. 
\subsection*{Scalar fields}
Let $f\in T^0 \SurfDomain$ be a scalar field with arbitrary smooth extension $f^e$, i.\,e.\ $f^e\vert_{\SurfDomain} = f$. The componentwise, tangential, or covariant derivative is given by:        \begin{align} \label{eq:scalar_grad}
    \GradC f = \GradP f = \GradSurf f
    &= (\ProjMat^e\nabla f^e)\vert_{\SurfDomain}
    = \ProjMat (\nabla f^e)\vert_{\SurfDomain} \ .
\end{align}
An alternative definition of the same operators without considering the extension $f^e$ is given by:
\begin{align}
    \GradC f = \GradP f = \GradSurf f
    &= g^{ij}\partial_j f \partial_i \param 
    = (g^{ij}\partial_j f \partial_i X^A) \eb_A \ ,
\end{align}
where $g^{ij}$ are the contravariant/inverse components of the metric tensor on $\SurfDomain$ obtained by a parametrization $\param$.

\subsection*{Vector fields}

Let $\vectvel = \vectvel_T + u_N \normalvec \in T\REAL^3\vert_{\SurfDomain}$ be a vector field with arbitrary smooth extension $\vectvel^e$, i.\,e.\ $\vectvel^e\vert_{\SurfDomain} = \vectvel$. We define the componentwise and tangential derivatives, respectively, by: 
\begin{align}
  \GradC\vectvel 
  &= ((\ProjMat^e\nabla) \vectvel^e)\vert_{\SurfDomain}
  = (\nabla \vectvel^e)\vert_{\SurfDomain}\ProjMat &\text{(componentwise derivative)}\,, \\
  \GradP\vectvel 
  &= (\ProjMat^e(\nabla \vectvel^e)\ProjMat^e)\vert_{\SurfDomain}
  = \ProjMat(\nabla \vectvel^e)\vert_{\SurfDomain}\ProjMat &\text{(tangential derivative)}.
\end{align}
A covariant derivative is uniquely defined for pure tangential vector fields (i.e., $u_N=0$) by: 
\begin{align}
    \GradSurf \vectvel_T 
    &= g^{jk}\left( \partial_k u_T^i + \Gamma_{kl}^i u_T^l \right) \partial_i \param \otimes \partial_j \param \ ,
\end{align}
where $\Gamma^{\cdot}_{\cdot\cdot}$ are the Christoffel symbols determined by the metric tensor. The relation to the tangential and componentwise derivatives is given by: 
\begin{align}
    \GradP\vectvel 
    &= g^{ik}g^{jl} \left( \partial_l \vectvel , \partial_k \param \right) \partial_i \param \otimes \partial_j \param
    = \ProjMat\GradC\vectvel
    = \GradSurf\vectvel_T - u_N \shapeOp\label{eq:vector_GradP}\,, \\
    \GradC\vectvel 
    &= \eb_A \otimes \GradSurf\vectvel^A
    = g^{ij}(\partial_j \vectvel) \otimes \partial_i \param
    = \GradSurf\vectvel_T - u_N \shapeOp + \normalvec\otimes\left( \GradSurf u_N + \shapeOp\vectvel_T \right) \label{eq:vector_GradC}\ .
\end{align}
By definition, the covariant and tangential derivatives of vector fields are purely tangential operators but the componentwise derivative contains normal components too. 
We define the corresponding divergence operators by:
\begin{align}
    \divS\vectvel_{T}
    &= \tr\GradSurf\vectvel_{T}
    =\partial_i u_T^i + \Gamma_{ik}^i u_T^k\,, \\
    \DivC\vectvel &= \DivP\vectvel 
    = \tr\GradC\vectvel = \tr\GradP\vectvel
    = \divS\vectvel_T - u_N\meanCurv \label{eq:vector_divCP} \formPeriod
\end{align}

\subsection*{Tensor fields}

Let $\StressTens \in T^2\REAL^3\vert_{\SurfDomain}$ be a tensor field with arbitrary smooth extension $\StressTens^e$, i.\,e.\ $\StressTens^e\vert_{\SurfDomain} = \StressTens$. The definition of the componentwise, tangential, or covariant derivative is possible for arbitrary tensor fields but we will show only the definitions of the operators used in this paper. To do so, we define the componentwise derivative by: 
\begin{align}
    \GradC\StressTens
    &= ((\ProjMat^e\nabla) \StressTens^e)\vert_{\SurfDomain}
    = (\nabla \StressTens^e)\vert_{\SurfDomain}\ProjMat\\
    &= \eb_A \otimes \eb_B \otimes \GradSurf\StressTensComp^{AB}
    = g^{ij}(\partial_j \StressTens) \otimes \partial_i \param \, .
\end{align}
The covariant divergence for the pure tangential tensor field $\boldsymbol{\epsilon} = \ProjMat\StressTens\ProjMat \in T^2 \SurfDomain$ is defined by:
\begin{align}
  \divS\boldsymbol{\epsilon}
  &=(\partial_j \epsilon^{ij} + \Gamma_{jk}^i \epsilon^{kj} + \Gamma_{jk}^j \epsilon^{ik}) \partial_i \param \ .
\end{align}
We define a componentwise divergence for right side tangential tensor field $\StressTens\ProjMat$ by: 
\begin{align}\label{eq:tensor_DivC}
    \DivC(\StressTens\ProjMat)
    &= \tr\GradC(\StressTens\ProjMat) 
    = \left( \eb_B , \GradSurf [\StressTens\ProjMat]^{AB} \right) \eb_A\\
    &= \divS( \ProjMat\StressTens\ProjMat) - \normalvec\StressTens\shapeOp
    + \left( ( \StressTens, \shapeOp ) + \divS( \normalvec\StressTens\ProjMat ) \right)\normalvec \notag\ .
\end{align}
The reason for defining only this restricted version of the componentwise divergence is that $\tr\circ\GradC$ would not be the adjoint operator of $\GradC$ and thus the integration by parts formula \eqref{eq:integration_by_parts_GradC_vector} would no longer hold in general.

\subsection*{Integration by parts formulae}
An integration by parts formulae can be derived for the differential operators above.
The following uses the global inner product neglecting the boundary terms.
For all tangential tensor fields $\boldsymbol{\epsilon} = \ProjMat\StressTens\ProjMat \in T^2 \SurfDomain$ and vectors fields $\vectvel  \in T\REAL^3\vert_{\SurfDomain} $, the following relations hold: 
\begin{align}
    \left( \StressTens, \GradP\vectvel \right) 
    = \left( \boldsymbol{\epsilon}, \GradP\vectvel \right)
    = \left( \boldsymbol{\epsilon}, \GradC\vectvel \right)
    &= - \left( \DivC \boldsymbol{\epsilon} , \vectvel \right) \label{eq:integration_by_parts_GradP_vector}\,, \\
   \left( \StressTens\ProjMat, \GradC\vectvel \right)
   &= - \left( \DivC (\StressTens\ProjMat) , \vectvel \right)    \label{eq:integration_by_parts_GradC_vector} \ .
\end{align}
These relations follow from \eqref{eq:vector_GradC} and \eqref{eq:tensor_DivC}.
For all $f\in T^0 \SurfDomain$ and $ \vectvel  \in T\REAL^3\vert_{\SurfDomain} $, we have:
\begin{align}
    \left( f , \DivP\vectvel \right) = \left( f , \DivC\vectvel \right)
    &= -\left( \DivC ( f\ProjMat ) , \vectvel \right) 
    = - \left( \GradSurf f + f\meanCurv\normalvec , \vectvel \right) \label{eq:integration_by_parts_divCP_vector}  \ ,
\end{align}
which follows from \eqref{eq:integration_by_parts_GradP_vector}, \eqref{eq:scalar_grad}, \eqref{eq:tensor_DivC},
and metric compatibility $\GradSurf\ProjMat=0$.

\subsection*{Laplace operators}
For scalar fields, we define a Laplace operator by the Laplace Beltrami operator $\LapSurf$:
\begin{align}
    \LapSurf f
    &= \divS\GradSurf f 
    = \DivP\GradP f
    = \DivC\GradC f
    =  g^{ij}\left( \partial_i\partial_j f - \Gamma_{ij}^{k}\partial_k f \right)\,,
\end{align}
and the componentwise Laplace operator for vector fields $\vectvel= u^A \eb_{A}$ by : 
\begin{align}
    \Delta_C \vectvel = \DivC \GradC \vectvel = (\LapSurf u^A)\eb_{A} \formComma\,,
\end{align}
where $u^A$ are the Cartesian components with respect to the Cartesian basis vectors $\eb_{A}$.

\section{Lagrange-D'Alembert principle for the full model}
\label{app:dalembert}

\subsection{General proceeding}
For state variables $\param$ (parametrization of $\SurfDomain$) and $\phi\in T^0\SurfDomain$ (phase field), 
and process variables  $\vectvel\in T\REAL^3\vert_{\SurfDomain}$ (material velocity) and $\dot{\phi}\in T^0\SurfDomain$ (phase field rate),
the Lagrange-D'Alembert principle generally reads (see \cite{marsden2013introduction}):
\begin{align}\label{eq:LagrangeDAlembert_gen_timeglobal}
    \left( \frac{\delta\mathcal{A}}{\delta\param}, \VecFieldY \right) + \left( \frac{\delta\mathcal{A}}{\delta\phi}, \psi \right)
    &= \int_{t_0}^{t_1} \left( \frac{\delta\DissPotential}{\delta\vectvel}, \VecFieldY \right) + \left( \frac{\delta\DissPotential}{\delta\dot{\phi}}, \psi \right)
    &\forall\VecFieldY\in T\REAL^3\vert_{\SurfDomain}, \psi\in T^0\SurfDomain\,,
\end{align}
where $ \DissPotential$ is the dissipation potential and 
$\mathcal{A} = \int_{t_0}^{t_1} \Lagrangian$ is the action functional, in which the 
Lagrangian is defined by $\Lagrangian= \energyKin - \energy$, with kinetic energy $\energyKin = \int_{\SurfDomain}\frac{\rho}{2}\NORM{\vectvel}^2$, potential energy $\energy$, and time interval $[t_0,t_1]$.

In a first step, we localize \eqref{eq:LagrangeDAlembert_gen_timeglobal} in time. 
Assuming that the potential energy $\energy$ depends only on state variables, we only have to consider the kinetic energy $\energyKin$ for this purpose.
Note that time integration commutes with spatial variations in a non-relativistic setting. 
For instance, the following holds 
$\left(\frac{\delta}{\delta\param}\int_{t_0}^{t_1} \energyKin , \VecFieldY \right) 
            = \int_{t_0}^{t_1}  \left( \frac{\delta\energyKin}{\delta\param} , \VecFieldY \right)$.
In contrast, spatial integration does clearly not commutes with spatial variation.
From \cite{Nitschke2022GoI}, we adopt the identity:
\begin{align}\label{eq:deformation_transport}
        \left(\frac{\delta }{\delta \param} \int_{\SurfDomain} f , \VecFieldY \right)
            &=\int_{\SurfDomain} \gradDef{\VecFieldY}f + f\DivP \VecFieldY 
        &\forall f \in T^0\SurfDomain, \VecFieldY\in T\REAL^3\vert_{\SurfDomain} \,,
\end{align}
where $\gradDef{\VecFieldY}:T^0\SurfDomain \rightarrow T^0\SurfDomain $ is the deformation derivative in direction of $\VecFieldY$, i.\,e.\ , 
\begin{align*}
    \gradDef{\VecFieldY}f
        &:= \frac{d}{d\varepsilon}\Big\vert_{\varepsilon=0} \left(f\vert_{\param+\varepsilon\VecFieldY} \right)
    &\forall f \in T^0\SurfDomain, \VecFieldY\in T\REAL^3\vert_{\SurfDomain} \,.    
\end{align*}
We assume a variational mass conservation, i.\,e.\ ,
\begin{align*}
    0=\left(\frac{\delta }{\delta \param} \int_{U} \rho , \VecFieldY \right)
        &=\int_{U} \gradDef{\VecFieldY}\rho + \rho\DivP \VecFieldY 
    &\forall U\subseteq\SurfDomain, \VecFieldY\in T\REAL^3\vert_{\SurfDomain} \,,
\end{align*}
valid for mass density $\rho\in T^0 \SurfDomain$.
Since the subset $U\subseteq\SurfDomain$ is arbitrary, $ \gradDef{\VecFieldY}\rho = -\rho\DivP \VecFieldY $ holds locally.
In addition, we assume temporal mass conservation, $\dot\rho = -\rho\DivP\vectvel$, and vanishing variation directions at times $t_0$ and $t_1$, 
i.\,e.\ $\VecFieldY\vert_{t_0} = \VecFieldY\vert_{t_1} = 0$.
The fundamental theorem of calculus, temporal integration by parts, transport formula, 
and $\gradDef{\VecFieldY}\vectvel = \dot{\VecFieldY} $ Cartesian-componentwise together yields:
\begin{align*}
    \left(\frac{\delta}{\delta\param}\int_{t_0}^{t_1} \energyKin , \VecFieldY \right)
        &= \int_{t_0}^{t_1}  \left( \frac{\delta\energyKin}{\delta\param} , \VecFieldY \right)
        = \int_{t_0}^{t_1} \left( \rho \vectvel, \dot{\VecFieldY} \right)
            + \int_{\SurfDomain} \frac{\NORM{\vectvel}^2}{2} \left( \gradDef{\VecFieldY}\rho + \rho\DivP\VecFieldY  \right)\\
        &= \int_{t_0}^{t_1} \frac{d}{dt} \left( \rho \vectvel, \VecFieldY \right)
            - \left( \rho \dot{\vectvel}, \VecFieldY \right) 
             -  \frac{1}{2}\int_{\SurfDomain} \left( \dot\rho + \rho\DivP\vectvel  \right)(\vectvel,\VecFieldY ) \\
        &= -\int_{t_0}^{t_1} \left( \rho \dot{\vectvel}, \VecFieldY \right)\,.
\end{align*}
Finally, the temporal localized version of the Lagrange-D'Alembert principle \eqref{eq:LagrangeDAlembert_gen_timeglobal} reads:
\begin{align}\label{eq:LagrangeDAlembert_gen}
    0
    &=  \left( \rho \dot{\vectvel} + \gradFunc{\param}\energy   + \gradFunc{\vectvel}\DissPotential  , \VecFieldY \right)
        + \left( \gradFunc{\phi}\energy   + \gradFunc{\dot{\phi}}\DissPotential, \psi \right)
    &\forall\VecFieldY\in T\REAL^3\vert_{\SurfDomain}, \psi\in T^0\SurfDomain \,,
\end{align}
where the negative generalized applied forces $ \gradFunc{\param}\energy, \gradFunc{\vectvel}\DissPotential \in T\REAL^3\vert_{\SurfDomain} $ 
and $\gradFunc{\phi}\energy, \gradFunc{\dot{\phi}}\DissPotential\in T^0\SurfDomain $ are given as $L^2$-gradients, 
i.\,e.\ $  ( \gradFunc{\param}\energy, \VecFieldY ) =  ( \frac{\delta\energy}{\delta\param}, \VecFieldY )  $
and $ ( \gradFunc{\vectvel}\DissPotential, \VecFieldY ) =  (\frac{\delta\DissPotential}{\delta\vectvel}, \VecFieldY ) $ 
for all virtual displacements $\VecFieldY\in T\REAL^3\vert_{\SurfDomain}$,
as well as $ ( \gradFunc{\phi}\energy, \psi ) = ( \frac{\delta\energy}{\delta\phi}, \psi ) $
and $ ( \gradFunc{\dot{\phi}}\DissPotential, \psi ) = ( \frac{\delta\DissPotential}{\delta\dot{\phi}}, \psi ) $
for all virtual displacements $ \psi\in T^0\SurfDomain $.

Next, in Section \ref{sec:var_of_pot}, we calculate the negative Lagrangian forces $ \gradFunc{\param}\energy$ and $ \gradFunc{\phi}\energy $, 
and the negative dissipative forces $ \gradFunc{\vectvel}\DissPotential $ and $ \gradFunc{\dot{\phi}}\DissPotential $ in Section \ref{sec:var_of_diss},
according to the potential energy $\energy$ and the dissipation potential $\DissPotential$ given in this paper.
Moreover, we implement local inextensibility of the material by the Lagrange-multiplier technique in Section \ref{sec:inext}.
Eventually, the strong formulation of \eqref{eq:LagrangeDAlembert_gen} and a tangential-normal splitting of the Lagrangian force  
to $\boldsymbol{b}_T := - \ProjMat \gradFunc{\param}\energy$ and  $\boldsymbol{b}_N := - (\normalvec, \gradFunc{\param}\energy)\normalvec$
leads to the full model given in Problem \ref{pb:fullmodel-dimensions}.
Note that, in our derivations, we assume a priori independence between the phase field $\phi$ and the surface given by the parametrization $\param$.
Therefore, $\gradDef{\VecFieldY}\phi = 0$ is valid for all $\VecFieldY\in T\REAL^3\vert_{\SurfDomain}$,
see \cite{Nitschke2022GoI} for more details, where this a priori condition is referred to the scalar gauge of surface independence.
In contrast to $L^2$-gradient flow techniques, here this assumption is not necessary to determine \eqref{eq:LagrangeDAlembert_gen},
since all terms containing $\gradDef{\VecFieldY}\phi$ would vanish anyway.
However, taken the generalized applied forces individually would comprises terms of the undetermined quantity $\gradDef{\VecFieldY}\phi$ in a weak sense.
Consequently, we would have to choose $\gradDef{\VecFieldY}\phi = 0$ to give this forces a determined meaning.

\subsection{Lagrangian forces}\label{sec:var_of_pot} 
In this section we consider the potential energy $\energy= \energyGL + \energyHelfr$, which is composed by the Ginzburg-Landau energy 
$\energyGL = \tcurve\int_{\SurfDomain} \frac{\interfaceparam}{2}\NORM{\GradSurf\CHsol}^2 +\frac{1}{\interfaceparam}\dWell(\CHsol)$ \eqref{eq:GL}
and the Helfrich energy 
$\energyHelfr = \int_{\SurfDomain} \frac{1}{2} \bendStiff(\CHsol)(\meanCurv-\meanCurv[0](\CHsol))^2$ \eqref{eq:helfrich}.
Without much mathematical effort, computing the negative generalized applied forces $\gradFunc{\phi}\energyGL, \gradFunc{\phi}\energyHelfr \in T^0\SurfDomain$ gives:
\begin{align*}
    \gradFunc{\phi}\energyGL
        &= \tcurve \left(-\interfaceparam\LapSurf\CHsol + \frac{1}{\interfaceparam}\dWell^{\prime}(\CHsol) \right)\,,\\
    \gradFunc{\phi}\energyHelfr
        &= \frac{1}{2}\bendStiff^{\prime}(\CHsol)(\meanCurv-\meanCurv[0](\CHsol))^2
                -\bendStiff(\CHsol)(\meanCurv-\meanCurv[0](\CHsol))\meanCurv[0]^{\prime}(\CHsol)\,.
\end{align*}

To derive the Ginzburg-Landau force $ -\gradFunc{\param}\energyGL $, we use that
$ \gradDef{\VecFieldY}g^{ij} = -( [\GradP\VecFieldY]^{ij} + [\GradP\VecFieldY]^{ji} ) $ is valid
for $g_{ij} = ( \partial_i\param,\partial_j\param )\ $ \cite{Nitschke2022GoI}.
Moreover, the relation $ \gradDef{\VecFieldY}\partial_i \phi = 0 $ holds,
according to the assumption $\gradDef{\VecFieldY}\phi = 0$.
This results in the following:
\begin{align*}
    \gradDef{\VecFieldY} \NORM{ \GradSurf\phi }^2
        &= \gradDef{\VecFieldY}\left( g^{ij} \partial_i \phi \partial_j \phi \right)
         = - \left(  \GradSurf\phi \otimes \GradSurf\phi , \GradP\VecFieldY + (\GradP\VecFieldY)^T \right)\\
        &=- 2 \left(  \GradSurf\phi \otimes \GradSurf\phi , \GradP\VecFieldY  \right)\,,
\end{align*}
by symmetry of the outer product.
With \eqref{eq:deformation_transport} and \eqref{eq:vector_divCP}, we obtain:
\begin{align*}
    \left( \frac{\delta\energyGL}{\delta\param}, \VecFieldY \right)
        &= \tcurve\int_{\SurfDomain} -\interfaceparam\left(  \GradSurf\phi \otimes \GradSurf\phi , \GradP\VecFieldY  \right) 
                 + \left(\frac{\interfaceparam}{2}\NORM{\GradSurf\CHsol}^2 +\frac{1}{\interfaceparam}\dWell(\CHsol)\right)\DivP \VecFieldY\\
         &= -\tcurve\left( \interfaceparam\GradSurf\phi \otimes \GradSurf\phi - \left( \frac{\interfaceparam}{2}\NORM{\GradSurf\CHsol}^2 +\frac{1}{\interfaceparam}\dWell(\CHsol)  \right)\ProjMat, \GradP\VecFieldY  \right)        \,.
\end{align*}
Eventually, integration by parts \eqref{eq:integration_by_parts_GradP_vector} yields:
\begin{align}
    \gradFunc{\param}\energyGL
        &= \tcurve\DivC \left( \interfaceparam\left(  \GradSurf\phi \otimes \GradSurf\phi - \frac{\NORM{\GradSurf\CHsol}^2}{2}\ProjMat \right) 
                        - \frac{1}{\interfaceparam}\dWell(\CHsol)\ProjMat\right) \in T\REAL^3\vert_{\SurfDomain} \,.
\end{align}
In terms of stress, this means that the Ginzburg-Landau energy induces a trace-free and symmetric tangential stress for phase separations 
and a volumetric stress for the double-well potential. 
To separate the induced  tangential and normal forces, we use metric compatibility $\GradSurf\ProjMat= 0$ 
and surface divergence \eqref{eq:tensor_DivC}.
This results in:
\begin{align}
    \ProjMat \gradFunc{\param}\energyGL
        &= \tcurve\divS \left( \interfaceparam\left(  \GradSurf\phi \otimes \GradSurf\phi - \frac{\NORM{\GradSurf\CHsol}^2}{2}\ProjMat \right) 
                        - \frac{1}{\interfaceparam}\dWell(\CHsol)\ProjMat \right) \notag\\
        &= \tcurve\left( \interfaceparam\LapSurf\CHsol - \frac{1}{\interfaceparam}\dWell^{\prime}(\CHsol) \right)\GradSurf\CHsol
         = -  (\gradFunc{\phi}\energyGL) \GradSurf\CHsol \in T \SurfDomain \,,\\
    \normalvec \gradFunc{\param}\energyGL
        &= \tcurve\left(  \interfaceparam\left(  \GradSurf\phi \otimes \GradSurf\phi - \frac{\NORM{\GradSurf\CHsol}^2}{2}\ProjMat \right) 
                        - \frac{1}{\interfaceparam}\dWell(\CHsol)\ProjMat , \shapeOp \right) \notag\\
        &= \tcurve\interfaceparam (\GradSurf\phi )\shapeOp \GradSurf\phi 
            - \tcurve\meanCurv\left( \frac{\interfaceparam}{2}\NORM{\GradSurf\CHsol}^2 +\frac{1}{\interfaceparam}\dWell(\CHsol) \right) \in T^0 \SurfDomain \,.
\end{align}

To calculate the Helfrich force $ -\gradFunc{\param}\energyHelfr $ we need the deformation derivative of the mean curvature.
From \cite{Nitschke2022GoI}, we already know the tangential part of the deformation derivative of the shape operator.
Explicitly, $ \ProjMat(\gradDef{\VecFieldY}\shapeOp)\ProjMat = \GradSurf ( \normalvec\GradC\VecFieldY ) - \shapeOp\GradP\VecFieldY$,
where $ \gradDef{\VecFieldY}\shapeOp $ is to be read Cartesian-componentwise.
Since $\tr$ and $\gradDef{\VecFieldY}$ commute in this way, 
we obtain $ \gradDef{\VecFieldY} \meanCurv = \divS( \normalvec\GradC\VecFieldY ) - (\shapeOp,\GradP\VecFieldY) $.
With \eqref{eq:deformation_transport} and integration by parts for the covariant divergence, we obtain:
\begin{align*}
    \left( \frac{\delta\energyHelfr}{\delta\param}, \VecFieldY \right)  
        &= \int_{\SurfDomain} \bendStiff(\CHsol)(\meanCurv-\meanCurv[0](\CHsol))
                \left( \gradDef{\VecFieldY} \meanCurv + \frac{1}{2}(\meanCurv-\meanCurv[0](\CHsol)) \DivP\VecFieldY \right) \\
        &= - \left( \normalvec\otimes\GradSurf( \bendStiff(\CHsol)(\meanCurv-\meanCurv[0](\CHsol)) ) 
                        +\bendStiff(\CHsol)(\meanCurv-\meanCurv[0](\CHsol)) \left( \shapeOp - \frac{\meanCurv-\meanCurv[0](\CHsol)}{2} \ProjMat \right) 
                , \GradC\VecFieldY\right)\,,
\end{align*}
and then:
\begin{align}
    \gradFunc{\param}\energyHelfr
        &= \DivC \left(  \bendStiff(\CHsol)(\meanCurv-\meanCurv[0](\CHsol)) \left( \shapeOp - \frac{\meanCurv-\meanCurv[0](\CHsol)}{2} \ProjMat \right) 
                        + \normalvec\otimes\GradSurf( \bendStiff(\CHsol)(\meanCurv-\meanCurv[0](\CHsol)) \right)\,,
\end{align}
by integration by parts \eqref{eq:integration_by_parts_GradC_vector}.
In terms of stress, we get an orthogonal decomposition of a trace-free and symmetric tangential stress tensor, 
a volumetric stress tensor scaled by $ \meanCurv[0](\CHsol) $,
and an additional stress tensor operating in the normal-tangential space $\normalvec\otimes T\SurfDomain$.
Using \eqref{eq:tensor_DivC}, 
metric compatibility $\GradSurf\ProjMat= 0$, the fact that $\shapeOp$ is curl-free and thus $\divS\shapeOp=\GradSurf\meanCurv$,
and $\GradSurf f(\CHsol) = f^{\prime}(\CHsol)\GradSurf\CHsol$ for $f\in\{ \bendStiff, \meanCurv[0] \}$, the
results in the tangential and normal Helfrich forces are:
\begin{align}
    \ProjMat \gradFunc{\param}\energyHelfr 
        &= \left(\bendStiff(\CHsol)(\meanCurv-\meanCurv[0](\CHsol))\meanCurv[0]^{\prime}(\CHsol)
            -\frac{1}{2}\bendStiff^{\prime}(\CHsol)(\meanCurv-\meanCurv[0](\CHsol))^2\right)
                \GradSurf\CHsol \notag\\
        &= - (\gradFunc{\phi}\energyHelfr) \GradSurf\CHsol  \in T \SurfDomain \,,\\
    \normalvec \gradFunc{\param}\energyHelfr 
        &= \LapSurf\left(\bendStiff(\CHsol)(\meanCurv-\meanCurv[0](\CHsol))\right)
         + \bendStiff(\CHsol)(\meanCurv-\meanCurv[0](\CHsol))(\NORM{\shapeOp}^2 - \frac{\meanCurv}{2}(\meanCurv-\meanCurv[0](\CHsol))) \,.
\end{align}

\subsection{Dissipative forces}\label{sec:var_of_diss}
In this section we consider the dissipation potential $\DissPotential = \DissPotential_V + \DissPotential_R + \DissPotential_\CHsol$,
where $\DissPotential_V = \int_{\SurfDomain} \frac{\viscosity}{4} \NORM{\GradP \vectvel + (\GradP \vectvel)^T}^2 $ 
is the viscous stress, $ \DissPotential_R = \int_{\SurfDomain} \frac{\KinCoef}{2}\NORM{\vectvel}^2 $ the friction with surrounding material,
and $\DissPotential_\CHsol =  \frac{1}{2\mobility} \left\vert \dot{\CHsol} \right\vert_{H^{-1}}^2 $ 
the ($H^{-1}$)-immobility potential.
Since the viscous stress potential does not depend on $\dot{\CHsol}$, we get:
\begin{align}
    \gradFunc{\dot{\phi}}\DissPotential_V &= 0 \,.
\end{align}
Variation with respect to $\vectvel$ gives
$ \left( \frac{\delta\DissPotential_V}{\delta\vectvel}, \VecFieldY \right)
        = \viscosity\left( \GradP \vectvel  + (\GradP \vectvel)^T , \GradP \VecFieldY\right) $
by symmetry.
Therefore, by applying integration by parts \eqref{eq:integration_by_parts_GradP_vector} we obtain the negative viscous force:
\begin{align}\label{eq:DuDissV}
    \gradFunc{\vectvel}\DissPotential_V
        &= - \viscosity\DivC \left( \GradP \vectvel  + (\GradP \vectvel)^T \right) \in T\REAL^3\vert_{\SurfDomain}\,,
\end{align}
containing twice the tangential strain rate tensor.
The friction potential does not depend on $\dot{\CHsol}$.
Therefore, we have:
\begin{align}
    \gradFunc{\dot{\phi}}\DissPotential_R 
        &= 0 \,,\\
    \gradFunc{\vectvel}\DissPotential_R \label{eq:DuDR}
        &= \KinCoef \vectvel \in T\REAL^3\vert_{\SurfDomain} \,.
\end{align}
We approach the immobility potential by defining a scalar field $\varphi[f]\in T^0 \SurfDomain$ implicitly so that it solve the equation
$\LapSurf \varphi[f] = f \in T^0 \SurfDomain$.
Due to this, we can write the immobility potential in a $L^2$-manner by 
$\DissPotential_\CHsol =  \frac{1}{2\mobility}  \int_{\SurfDomain}  \NORM{\GradSurf \varphi[\dot{\CHsol}]}^2   $.
Since $\varphi[f]$ is linear in $f$, this yields:
\begin{align*}
    \left(  \frac{\delta\DissPotential_\CHsol}{\delta\dot{\phi}}, \psi \right)
        &= \frac{1}{\mobility} \left( \GradSurf\varphi[\dot{\CHsol}] , \GradSurf\varphi[\psi] \right)
         = - \frac{1}{\mobility} \left( \varphi[\dot{\CHsol}] , \LapSurf\varphi[\psi] \right)
         = - \frac{1}{\mobility} \left( \varphi[\dot{\CHsol}] , \psi \right) \,.
\end{align*}
By considering $\gradDef{\VecFieldY}\phi = 0$, we get the dissipative forces: 
\begin{align}
    \gradFunc{\dot{\phi}}\DissPotential_\CHsol 
        &= - \frac{1}{\mobility}\varphi[\dot{\CHsol}]\,,  
    &(\mbox{with }\ \LapSurf \varphi[\dot{\CHsol}] &= \dot{\CHsol})  \,,\\
    \gradFunc{\vectvel}\DissPotential_\CHsol
        &= 0 \,.
\end{align}
Note that for all sufficient smooth $f \in T^0 \SurfDomain$ and $\mobility\neq0$ the equation $ \gradFunc{\dot{\phi}}\DissPotential_\CHsol + f =0  $ implies
$ \dot{\CHsol} = \mobility \LapSurf f  $.

\subsection{Implementation of local inextensibility}\label{sec:inext}
The Lagrange-D'Alembert principle \eqref{eq:LagrangeDAlembert_gen} holds by assuming mass conservation $ 0 = \dot\rho + \rho\DivP\vectvel$.
If we also want to set $\dot\rho=0$, or equivalently $\DivP\vectvel=0$, we have to constrain the solution space of \eqref{eq:LagrangeDAlembert_gen}.
We obtain this constraint by the Lagrange multiplier technique, where the Lagrange function 
$ \mathcal{C}_{IE}=-\int_{\SurfDomain}\press\DivP\vectvel$
has to be included appropriately into the variation process. 
The Lagrange multiplier $\press\in T^0 \SurfDomain$ yields a new degree of freedom, which we associate formally to a process variable,
since its purpose is to constrain the process variable $\vectvel$.
Therefore, we obtain additional negative applied forces:
\begin{align}
    \gradFunc{\press} \mathcal{C}_{IE} 
        &= - \DivP\vectvel \in T^0 \SurfDomain\,,\\
    \gradFunc{\vectvel} \mathcal{C}_{IE} 
        &= \DivC( \press \ProjMat ) 
        = \GradSurf \press +  \press\meanCurv\normalvec  \in T\REAL^3\vert_{\SurfDomain} \, ,
\end{align}
by 
$ ( \gradFunc{\press} \mathcal{C}_{IE}, q ) :=  (\frac{\delta\mathcal{C}_{IE}}{\delta\press}, q ) $ and 
$ ( \gradFunc{\vectvel} \mathcal{C}_{IE}, \VecFieldY ) :=  (\frac{\delta\mathcal{C}_{IE}}{\delta\vectvel}, \VecFieldY ) $,
with virtual displacements $ q \in T^0 \SurfDomain $ and $ \VecFieldY \in T\REAL^3\vert_{\SurfDomain} $,
and integration by parts \eqref{eq:integration_by_parts_divCP_vector}.
We do not consider a variation with respect to $\dot{\CHsol}$, since this would result in $ \gradFunc{\dot{\CHsol}} \mathcal{C}_{IE} = 0$ anyway.
Including this into \eqref{eq:LagrangeDAlembert_gen} gives:
\begin{align*}
    0
    &=  \left( \rho \dot{\vectvel} + \gradFunc{\param}\energy   + \gradFunc{\vectvel}\DissPotential  +  \gradFunc{\vectvel} \mathcal{C}_{IE}, \VecFieldY \right)
        + \left( \gradFunc{\phi}\energy   + \gradFunc{\dot{\phi}}\DissPotential, \psi \right)
        + \left( \gradFunc{\press} \mathcal{C}_{IE}, q \right) \,,
\end{align*}
for all $ \VecFieldY\in T\REAL^3\vert_{\SurfDomain}$ and $ \psi,q \in T^0\SurfDomain $.
Mutual independence of the virtual displacements leads to the strong formulation:
\begin{align}
    \rho\dot{\vectvel} \label{eq:app_NSE}
        &= -\left( \DivC( \press \ProjMat ) + \gradFunc{\param}\energyGL + \gradFunc{\param}\energyHelfr + \gradFunc{\vectvel}\DissPotential_{V} + \gradFunc{\vectvel}\DissPotential_{R} \right)\,, \\
    \dot{\CHsol}
        &= \mobility \LapSurf \left( \gradFunc{\CHsol} \energyGL +  \gradFunc{\CHsol} \energyHelfr \right)\,, \\
    0 &= \DivP\vectvel \label{eq:app_MC} \,,
\end{align}
recalling the calculations in Sections \ref{sec:var_of_pot} and \ref{sec:var_of_diss}.
Since \eqref{eq:app_NSE} and \eqref{eq:app_MC} constitute the inextensible Navier-Stokes equations with some additional forces,
it is justified to call the Lagrange multiplier $\press$ the pressure.

\subsection{Total energy rate}\label{sec:energy_rate}
In this section we investigate the total energy rate 
$\frac{d}{dt} ( \energyKin + \energy ) $ comprising the sum of the kinetic energy $\energyKin$ and the potential energy $\energy= \energyGL + \energyHelfr$.
For eqs.  \eqref{eq:app_NSE} -- \eqref{eq:app_MC} to be thermodynamically consistent,
it is necessary that the total energy rate only depends on the dissipation potential $\DissPotential = \DissPotential_V + \DissPotential_R + \DissPotential_\CHsol$,
since this is the only mechanism that allows energy exchange with the surrounding of $\SurfDomain$,
and the total energy $\energyKin + \energy$ must not increase with time.

Using \eqref{eq:app_NSE} -- \eqref{eq:app_MC}, including $\dot{\rho}=0$, and applying the chain rule yields:
\begin{align*}
    \frac{d}{dt} \left( \energyKin + \energy \right)
        &= \left( \rho \dot{\vectvel}  + \gradFunc{\param}\energy , \vectvel \right)
            + \left( \gradFunc{\CHsol}\energy, \dot{\CHsol} \right)\\
        &= - \left( \DivC( \press \ProjMat ) + \gradFunc{\vectvel} \DissPotential_V +  \gradFunc{\vectvel} \DissPotential_R , \vectvel \right)
           - \frac{1}{M} \left\vert \dot{\CHsol} \right\vert_{H^{-1}}^2 \,.
\end{align*}
Integration by parts \eqref{eq:integration_by_parts_divCP_vector} and \eqref{eq:integration_by_parts_GradC_vector},  
$ \gradFunc{\vectvel} \DissPotential_V $ \eqref{eq:DuDissV}, and $\gradFunc{\vectvel} \DissPotential_R$ \eqref{eq:DuDR} results in:
\begin{align*}
    \left( \DivC( \press \ProjMat ) , \vectvel \right)
        &= -\left( \press, \DivP\vectvel \right) = 0 \,,\\
    \left(  \gradFunc{\vectvel} \DissPotential_V , \vectvel \right)
        &= \viscosity\left( \GradP\vectvel + (\GradP\vectvel)^T  , \GradP\vectvel \right)
         = \frac{\viscosity}{2} \int_{\SurfDomain} \NORM{\GradP\vectvel + (\GradP\vectvel)^T}^2 \,, \\
    \left(  \gradFunc{\vectvel} \DissPotential_R , \vectvel \right)
        &= \KinCoef\int_{\SurfDomain} \NORM{\vectvel}^2 \,.
\end{align*}
As a consequence, we finally obtain:
\begin{align}\label{eq:energy_rate}
    \frac{d}{dt} \left( \energyKin + \energy \right)
        &= -2 \DissPotential \leq 0 \,,
\end{align}
which satisfies our requirements above.

\subsection{Material time derivative}\label{sec:material_derivative}

In this section we explain the material time derivative in details.
Following the description in \cite{nitschke2022TimeDerivative} for scalar fields $\phi\in T^0\SurfDomain$, 
the material derivative can be written as:
\begin{align}\label{eq:material_derivative}
    \dot{\phi} &= \partial_t \phi + \nabla_{\vectvelV}\phi \,,
\end{align}
where $\vectvelV\in T^1 \SurfDomain$ is the so-called relative velocity.
To be able to evaluate both summands on the right-hand-side, we use the two different parametrizations
$\param:(t,y^1,y^2) \mapsto \param(t,y^1,y^2)\in\SurfDomain(t)$ and
$\param_{\mfrak}:(t,y_{\mfrak}^1,y_{\mfrak}^2) \mapsto \param_{\mfrak}(t,y_{\mfrak}^1,y_{\mfrak}^2)\in\SurfDomain(t)$.
Both parametrizations describe the same moving surface, where $\param_{\mfrak}$ describe the material,
\ie\ the map $\param_{\mfrak}\vert_{(y_{\mfrak}^1,y_{\mfrak}^2)}:[t_0,t_1] \rightarrow \REAL^3 $ details the path of a single material particle in time.
In contrast, $\param$ is arbitrary as long as it depicts the same moving surface.
This means that $\param$ acts as an observer of $\SurfDomain$.
With this clarified, the scalar field $\phi$ can be represented by evaluating $\phi(t,y^1,y^2)$ as well as $\tilde{\phi}(t,\xb(t))$,
where $\xb(t)=\param(t,y^1,y^2)$ or $(y^1,y^2)=\param\vert_{t}^{-1}(\xb(t))$, respectively. 
To make it clear, $\phi$ and $\tilde{\phi}$ are the very same scalar field from a physical point of view, only the mathematical representation differs. 
Partial derivatives depend on the mathematical description of their argument, contrarily to total derivatives, for example.
Hence, we have to be careful when evaluating the first summand in \eqref{eq:material_derivative}. 
This partial time derivative reads:
\begin{align*}
    \partial_t \phi( t,y^1,y^2)  
        &= \frac{d}{d\tau}\Big\vert_{\tau=0} \phi( t+\tau,y^1,y^2)
         = \frac{d}{d\tau}\Big\vert_{\tau=0} \tilde{\phi}( t+\tau, \param(t+\tau,y^1,y^2)) \formComma
\end{align*}
or in terms of a differential quotient:
\begin{align*}
    \widetilde{\partial_t \phi}( t, \xb(t)) 
        &= ( \partial_t \phi( t,y^1,y^2) )\vert_{(y^1,y^2)=\param\vert_{t}^{-1}(\xb(t))}\\
        &= \frac{1}{\tau} \left(  \tilde{\phi}( t+\tau, \xb(t+\tau)) - \tilde{\phi}( t, \xb(t)) \right) +\mathcal{O}(\tau)\,,
\end{align*}
with respect to global coordinates $\xb(t)\in\SurfDomain(t)$, where $\widetilde{\partial_t \phi}$ is equal to $\partial_t \phi$ due to the relation $\xb(t)=\param(t,y^1,y^2)$.
As a consequence, $\partial_t \phi$ is well defined and can be discretize in time with respect to global coordinates.
Since $\partial_t \phi$ represents only the observer rate of $\phi$, the second summand in \eqref{eq:material_derivative} represents the correction to obtain the material rate.
The relative velocity is given by $\wb= \ub - \vb$, where $\ub$ is the material velocity and $\vb$ the observer velocity.
In local observer coordinates, they are given by:
\begin{align*}
    \ub(t,y^1,y^2)
        &= \partial_t\param_{\mfrak}( t, \param_{\mfrak}\vert_{t}^{-1}(\param (t,y^1,y^2)) ) \formComma
    &\vb(t,y^1,y^2)
        &= \partial_t \param (t,y^1,y^2) \formPeriod
\end{align*}
As a consequence, we have:
\begin{align*}
    \tilde{\ub}(t, \xb(t)) 
        &= \ub(t,y^1,y^2)\vert_{(y^1,y^2)=\param\vert_{t}^{-1}(\xb(t))}
         = ( \partial_t \xb_{\mfrak}(t) )\vert_{\xb_{\mfrak}(t) = \xb(t)} \formComma\\
    \tilde{\vb}(t, \xb(t)) 
        &= \vb(t,y^1,y^2)\vert_{(y^1,y^2)=\param\vert_{t}^{-1}(\xb(t))}
         = \partial_t \xb(t) \formComma
\end{align*}
where the relation between global and material local coordinates is given by $\xb_{\mfrak}(t)=\param_{\mfrak}(t,y_{\mfrak}^1,y_{\mfrak}^2)$ or $(y_{\mfrak}^1,y_{\mfrak}^2)=\param_{\mfrak}\vert_{t}^{-1}(\xb_{\mfrak}(t))$, respectively.
Note that for geometrical reasons, $(\vb,\normalvec)=(\ub,\normalvec)$ holds.
Due to this, the relative velocity is a tangential vector field, \ie\ $\wb=\ub-\vb\in T^1 \SurfDomain$,
and therefore $\nabla_{\wb}\phi = ( \GradSurf\phi, \wb )$ is valid.
Putting all together, we can write the material derivative \eqref{eq:material_derivative} 
for  $\tilde{\dot{\phi}}(t, \xb(t))  = \dot{\phi}(t,y^1,y^2)\vert_{(y^1,y^2)=\param\vert_{t}^{-1}(\xb(t))}$ as:
\begin{align}\label{eq:material_derivative_global_diffquot}
    \tilde{\dot{\phi}}(t, \xb(t)) 
        &= \frac{1}{\tau} \left(  \tilde{\phi}( t+\tau, \xb(t+\tau)) - \tilde{\phi}( t, \xb(t)) \right) \notag\\
        &\quad +\left(  \GradSurf\tilde{\phi}(t, \xb(t)) \ , \ \tilde{\ub}(t, \xb(t)) 
                                            - \frac{1}{\tau} \left( \xb(t+\tau) - \xb(t) \right) \right)
            +\mathcal{O}(\tau) 
\end{align}
in terms of global coordinates and differential quotients.
This representation is very suitable to implement flow equations on a surface grid that does not necessarily have to follow the material flow. 
Here, $\xb(t)$ is given by the surface grid coordinates and the material velocity $\tilde{\ub}$ is part of the solution defined on these grid coordinates. To simplify the notation in all the other sections, we do not make use of the tilde, since all the field quantities with or without the tilde are exactly the same fields. They differ only within the argument formulation, which is not used in any case outside of this section.
Note that we can also apply \eqref{eq:material_derivative}, or \eqref{eq:material_derivative_global_diffquot} respectively, for the material acceleration $\dot{\ub}$ with respect to a Cartesian frame ${\eb_I}$.
Since $\dot{\eb}_I=0$, the relation $\dot{\ub} = \dot{u}^I \eb_I$ holds, for $ \ub = u^I \eb_I$.

There are other representations of the material derivative \eqref{eq:material_derivative}.
For instance \cite{art:Dziuk2013} stated:
\begin{align}\label{eq:material_derivative_dziuk}
    \hat{\dot{\phi}}(t,\xb_{\mfrak}(t)) 
        &= \partial_t \hat{\phi} (t,\xb_{\mfrak}(t))  + \left( \nabla \hat{\phi}, \hat{\ub}(t,\xb_{\mfrak}(t))  \right)\,,
\end{align}
when re-written in our notation\footnote{In \cite{art:Dziuk2013}, the authors use a function $G(t,\cdot): \SurfDomain(0) \rightarrow \SurfDomain(t)$. But we can define this map by $G(t,\xb_{\mfrak}(0)) := \xb_{\mfrak}(t)$ implicitly.}, where 
$ \hat{\phi}(t, \xb_{\mfrak}(t) )  = \phi(t, y^1,y^2) $ 
for $ (y^1,y^2)=\param\vert_t^{-1}(\xb_{\mfrak}(t)) $ and 
$ \hat{\ub}(t,\xb_{\mfrak}(t)) = \partial_t\xb_{\mfrak}(t)  $.
As already mentioned in \cite{art:Dziuk2013},  both summands on the right-hand-side are not defined on a geometrically nonstationary surface, if considered individually. 
For the partial time derivative we have to evaluate $\hat{\phi} (t+\tau,\xb_{\mfrak}(t))$ for a small time step $\tau$.
However, $\hat{\phi}\vert_{t+\tau}$ is only defined at future locations $\xb_{\mfrak}(t+\tau)\in\SurfDomain(t+\tau)$ and not at the current $\SurfDomain(t)$ if no further assumptions are considered.
The same applies to the second summand, which depends on the normal derivative $ ( \nabla \hat{\phi}, \normalvec ) $ if the material velocity $\hat{\ub}$ comprises a normal velocity.
Such a normal derivative cannot be obtained without considering its extensions outside of the surface.
One might be tempted to interpret \eqref{eq:material_derivative_dziuk} as a fully Eulerian perspective, since 
$(\partial_t \hat{\phi})\vert_{\xb_{\mfrak}(t)=\param(t,y^1,y^2)} = \partial_t \phi$ is valid for $\partial_t\param = 0$.
However, a fully Eulerian surface observer does not exists if the surface comprises a flow in normal direction.
Any normal part of the surface motion has to be treated in a Lagrangian perspective at least.
In general, \eqref{eq:material_derivative_dziuk} is equal to \eqref{eq:material_derivative},
since $(\partial_t \hat{\phi})\vert_{\xb_{\mfrak}(t)=\param(t,y^1,y^2)} = \partial_t \phi - ( \nabla\phi , \vb ) $ holds formally in arbitrary observer coordinates.

\section{Non-dimensionalization} \label{app:dedim}

In order to non-dimensionalize the model in Problem \ref{pb:fullmodel-dimensions}, we write the rescaled variables by using the symbol $\hat{\cdot}$. We define the rescaled parametrization $\param = \hat{\param} L$, the velocity $\vectvel = \hat{\vectvel}U$, the mean curvature $\meanCurv = {\hat{\meanCurv}}/{L}$, the phase field $\phi=\hat\phi$, and the chemical potential $\chempot = \rho\hat\chempot U^2$ for the characteristic length $L$ and a characteristic velocity $U$, which defines a rescaled time $t={L}/{U}\hat{t}$. According to the parametrization, we rescale the gradient of the embedded space by $\nabla = {\hat{\nabla}}/{L}$. This rescales the surface operators accordingly, e.g., the material time derivative reads $(\partial_t + \nabla_{\vectvelV})= {U}/{L}(\partial_{\hat t}+\hat{\nabla}_{\hat\vectvelV})$. We introduce the rescaled material functions: $\hat{\bendStiff}(\CHsol) = \bendStiff(\CHsol) / (\rho U^2 L^2$), $\hat{\meanCurv}_0(\CHsol) = {\meanCurv}_0(\CHsol)L$, and $\hat{W}(\CHsol) = {1}/{L^2}W(\CHsol)$; and the rescaled material parameters: $\Reynolds = {\rho LU}/{\viscosity}$, denoting the Reynolds number, $\hat{m} = {\rho m U}/{L}$, $\hat\sigma = \tilde\sigma/(\rho U^2L^2)$, and $\hat\gamma = {\gamma L}/{(\rho U)}$. Considering the rescaled variables and replacing the material functions and parameters into the model, we obtain the following non-dimensional system:
\begin{align*}
  \partial_{\hat{t}}{\hat\CHsol}+ \hat\nabla_{\hat\vectvelV}\CHsol &= \hat m\,\hat\LapSurf\hat\chempot \,,\\
  \hat\chempot&=-\hat\sigma\interfaceparam\,\hat\LapSurf\CHsol+\frac{\tcurve}{\interfaceparam}\hat\dWell^{\prime}(\hat\CHsol)
  +\frac{\hat\bendStiff^{\prime}}{2}(\hat\CHsol)\left(\hat\meanCurv-\hat\meanCurv_0(\CHsol)\right)^2
  -\hat\bendStiff(\hat\CHsol)\hat\meanCurv^{\prime}_0(\hat\CHsol)\left(\hat\meanCurv-\hat\meanCurv_0(\hat\CHsol)\right)\,,\\
  \partial_{\hat{t}} \hat\vectvel + \hat\nabla_{\hat\vectvelV}\hat\vectvel &= -\hat\GradSurf p -
  p\hat\meanCurv\normalvec + \frac{2}{\Reynolds}\hat\DivC\hat\StressTens - \hat\gamma \hat\vectvel + \hat{\boldsymbol{b}}_{T}+
  \hat{\boldsymbol{b}}_N \,,\\
  \hat{\DivP}\hat\vectvel &= 0 \,,
\end{align*}
where $\hat{\boldsymbol{b}}_{T}$ and $\hat{\boldsymbol{b}}_{N}$ correspond to the rescaled tangential and normal bending terms, respectively. 
In the main sections of the paper, we drop the symbol $\hat{\cdot}$ for a better readability.

\section{Derivation of the overdamped limit model}\label{app:overdamp}
We consider the overdamped limit of the model in Problem \ref{pb:fullmodel}. We define the rescaled parameters $\delta = 1/{\gamma}$, $\tilde{\tcurve} = \delta\tcurve$,  $\tilde{\bendStiff} = \delta\bendStiff$, and $\tilde{\mobility} ={\mobility}/{\delta}$, and the rescaled variables $\tilde{p} = \delta p$ and $\tilde{\chempot}= \delta\chempot$. Substituting these parameters and variables into the model, we get:
\begin{align*}
  \partial_t{\CHsol}+ \nabla_{\vectvelV}\CHsol =& \tilde\mobility\,\LapSurf\tilde\chempot \,,\\
  \tilde\chempot=&\tilde\tcurve \left(-\interfaceparam\,\LapSurf\CHsol+\frac{1}{\interfaceparam}\dWell^{\prime}(\CHsol)\right)
  +\frac{1}{2}\tilde\bendStiff^{\prime}(\CHsol)\left(\meanCurv-\meanCurv[0](\CHsol)\right)^2 \\
  &-\tilde\bendStiff(\CHsol)\meanCurv^{\prime}_0(\CHsol)\left(\meanCurv-\meanCurv[0](\CHsol)\right)\,, \\
  \delta\partial_t \vectvel + \delta\GradSurfConv{\vectvelV}\vectvel =& -\GradSurf \tilde{p}
  -\tilde{p}\meanCurv\normalvec + \frac{2\delta}{\Reynolds}\DivP\StressTens - \vectvel + \tilde{\boldsymbol{b}}_{T}+
  \tilde{\boldsymbol{b}}_N \,,\\
  \DivP\vectvel =& 0 \,,
\end{align*}
where $\tilde{\boldsymbol{b}}_{T}$ and $\tilde{\boldsymbol{b}}_N$ are the tangential and normal forces with respect to the rescaled variables and $\vectvelV=\vectvel-\partial_t\param$ as before. The overdamped limit is then obtained by considering $\delta \searrow 0$, and it reads:
\begin{align}
  \partial_t{\CHsol}+ \nabla_{\vectvelV}\CHsol =& \tilde\mobility\,\LapSurf\tilde\chempot \,,\\
  \tilde\chempot=& \tilde\tcurve \left(-\interfaceparam\,\LapSurf\CHsol+\frac{1}{\interfaceparam}\dWell^{\prime}(\CHsol)\right)
  +\frac{1}{2}\tilde\bendStiff^{\prime}(\CHsol)\left(\meanCurv-\meanCurv[0](\CHsol)\right)^2 \nonumber\\
  &-\tilde\bendStiff(\CHsol)\meanCurv^{\prime}_0(\CHsol)\left(\meanCurv-\meanCurv[0](\CHsol)\right) \,,\\
  \vectvel =& -\GradSurf \tilde{p} -\tilde{p}\meanCurv\normalvec + \tilde{\boldsymbol{b}}_T +\tilde{\boldsymbol{b}}_N \label{subeq:vel} \,,\\
  \DivP\vectvel =& 0 \label{subeq:conti}\,.
\end{align}
Considering $\vectvel_T = \ProjMat\vectvel$ and $u_N = \vectvel\cdot\normalvec$, eq. \eqref{subeq:vel} can be split into a tangential and a normal part. Substituting these terms in eq. \eqref{subeq:conti} and using the definitions for 
$\tilde{\boldsymbol{b}}_T$ and $ \tilde{\boldsymbol{b}}_N$ we obtain: 
\begin{align}
  \begin{aligned} \label{app:eq:chbending}
  \partial_t{\CHsol}+ \nabla_{\vectvelV}\CHsol =& \tilde\mobility\,\LapSurf\tilde\chempot \,,\\
  \tilde\chempot=&\tilde\tcurve \left(-\interfaceparam\,\LapSurf\CHsol+\frac{1}{\interfaceparam}\dWell^{\prime}(\CHsol)\right)
  +\frac{1}{2}\tilde\bendStiff^{\prime}(\CHsol)\left(\meanCurv-\meanCurv[0](\CHsol)\right)^2 \\
  &-\tilde\bendStiff(\CHsol)\meanCurv^{\prime}_0(\CHsol)\left(\meanCurv-\meanCurv[0](\CHsol)\right) \,,\\
  u_N =& -\tilde{p}\meanCurv +\tilde{b}_N \,,\\
  \vectvel_T =& -\GradSurf \tilde{p} + \tilde\chempot\GradSurf\CHsol \,,\\
  -\LapSurf\tilde{p} + \tilde{p}\meanCurv^2 =& -\divS \left( \tilde{\chempot}\GradSurf\CHsol\right) + \tilde{b}_N\meanCurv\,,
  \end{aligned}
\end{align}
where $\tilde{b}_N$ is given by:
\begin{align*}
\tilde{b}_N =& \tilde{\boldsymbol{b}}_N\cdot\normalvec \\
=& -\LapSurf(\bendStiff(\CHsol)(\meanCurv-\meanCurv[0](\CHsol)))
  - \bendStiff(\CHsol)(\meanCurv-\meanCurv[0](\CHsol))\left(\Vert\shapeOp\Vert^2-\frac{1}{2}\meanCurv(\meanCurv-\meanCurv[0](\CHsol)) \right)
  \\ 
  &+\tilde{\tcurve}\left(\frac{\interfaceparam}{2}\Vert\GradSurf\CHsol\Vert^2+\frac{1}{\interfaceparam}
  W(\CHsol)\right)\meanCurv   -\tilde{\tcurve}\interfaceparam\GradSurf\CHsol^T\shapeOp\GradSurf\CHsol \,.
\end{align*}
This model is closely related to a model discussed in~\cite{Hauser2013}. To perform the comparison, we recall that the total time derivative of the potential energy considered in \ref{sec:energy_rate} is given by:
\begin{equation}\label{app:eq:transport}
    \frac{d}{dt} \energy = (\dot{\phi},\gradFunc{\phi}\energy) + (\vectvel,\gradFunc{\param}\energy).
\end{equation}
In \cite{Hauser2013}, the material derivative $\dot{\phi}$ is reduced to the partial time derivative $\partial_t\phi$, and this leads to the differences noted between the models.

More commonly used models reduce the local inextensibility constraint and only consider a global area constraint: $$\int_{\SurfDomain} \DivP\vectvel \mathrm{d}\SurfDomain= \int_{\SurfDomain} \divS\vectvel_T - u_N\meanCurv\mathrm{d}\SurfDomain = -\int_{\SurfDomain} u_N\meanCurv\mathrm{d}\SurfDomain = 0\,,$$ which can be realized by: 
\begin{align*}
  \tilde{p} = \frac{\int_{\SurfDomain} u_N\meanCurv\mathrm{d}\SurfDomain}{\int_{\SurfDomain}\meanCurv\mathrm{d}\SurfDomain} \,.
\end{align*}
In this case $\tilde{p}$ and $\vectvel_T$ are independent. The corresponding system reads:
\begin{align*}
  \begin{aligned} 
  \partial_t{\CHsol} + \nabla_{\vectvelV}\CHsol =& \tilde\mobility\,\LapSurf\tilde\chempot \,,\\
  \tilde\chempot=&\tilde\tcurve\left(-\interfaceparam\,\LapSurf\CHsol+\frac{1}{\interfaceparam}\dWell^{\prime}(\CHsol)\right)
  +\frac{1}{2}\tilde\bendStiff^{\prime}(\CHsol)\left(\meanCurv-\meanCurv[0](\CHsol)\right)^2 \\
  &-\tilde\bendStiff(\CHsol)\meanCurv^{\prime}_0(\CHsol)\left(\meanCurv-\meanCurv[0](\CHsol)\right) \,,\\
  u_N =& -\tilde{p}\meanCurv +\tilde{b}_N\,,\\
  \vectvel_T =& \tilde\chempot\GradSurf\CHsol \,,\\
  \tilde{p} =& \frac{\int_{\SurfDomain} \tilde{b}_N\meanCurv\mathrm{d}\SurfDomain}{\int_{\SurfDomain}\meanCurv\mathrm{d}\SurfDomain} \, ,
  \end{aligned}
\end{align*}
with $\vectvelV = \vectvel_T + u_N\normalvec-\partial_t\param$. This system has also been considered in~\cite{Hauser2013} and shows the same conceptional differences discussed above. Other approaches, with $\meanCurv[0] = 0$, $\tilde\bendStiff(\CHsol) = \tilde\bendStiff$, and that use a penalization approach to enforce the global area constraint instead of the Lagrange multiplier, have been considered in \cite{elliot_stinner2009}.

\bibliography{strings,biblio}
\bibliographystyle{plain}

\end{document}